\documentclass[12pt]{amsart}
\usepackage{amsmath,amssymb,latexsym,psfig}

\newcommand{\bC}{\mbox{${\mathbb C}$}}

\newcommand{\br}[1]{\langle#1\rangle}

\newcommand{\dual}[2]{\br{#1 \:|\: #2}}

\newcommand{\setp}[2]{\mbox{$\{#1\::\:#2\}$}}
\newcommand{\casec}[5]{#1#2 \left\{ \begin{array}{ll} #3 &\mbox{if}\spa #4 \\ #5 &\mbox{otherwise\,,} \end{array} \right.}
\newcommand{\casesc}[5]{#1#2 \left\{ \begin{array}{ll} #3 &\mbox{if}\spa #4 \\ #5 &\mbox{otherwise\,;} \end{array} \right.}
\newcommand{\stacksum}[2]{\sum_{\begin{array}{c}\vspace{-7mm}\;\\ \vspace{-2mm}\scriptstyle{#1}\\ \scriptstyle{#2}\end{array}} }

\newcommand{\CP}{\bC P^\infty}

\newcommand{\calp}{{\cal P}}
\newcommand{\cald}{M}
\newcommand{\cals}{{\cal S}}
\newcommand{\calw}{{\cal W}}
\newcommand{\calh}{{\cal H}}
\newcommand{\bF}{{\mathbb F}}

\newcommand{\nub}{\overline{\nu}}

\newcommand{\hgt}{\operatorname{ht}}

\makeatletter
\let\choose\@@choose
\let\cal\mathcal
\makeatother

\makeatletter
\let\cal\mathcal
\makeatother

\newtheorem{conj}[equation]{Conjecture}
\newtheorem{prop}[equation]{Proposition}

\newtheorem{thm}[equation]{Theorem}
\newtheorem{rem}[equation]{Remarks}
\newtheorem{rema}[equation]{Remark}
\newtheorem{cor}[equation]{Corollary}
\newtheorem{lem}[equation]{Lemma}
\newtheorem{exa}[equation]{Example}

\makeatletter\@addtoreset{equation}{section}\makeatother

\newcommand{\spa}{\;\: }
\newcommand{\bQ}{{\mathbb Q}} 
\newcommand{\bZ}{{\mathbb Z}}
\newcommand{\case}[5]{#1#2 \left\{ \begin{array}{ll} #3 &\mbox{if $#4$} \\#5 &\mbox{otherwise\,.} \end{array} \right.}
\newcommand{\casel}[5]{#1#2 \left\{ \begin{array}{ll} #3 &\mbox{if $#4$} \\ \;\\#5 &\mbox{otherwise\,;} \end{array} \right.}
\newcommand{\caset}[7]{#1#2 \left\{ \begin{array}{lll} #3 &\mbox{if $#4$} \\ #5&\mbox{if $#6$}\\#7 &\mbox{otherwise\,;} \end{array} \right.}
\newcommand{\map}[3]{\mbox{$#1 \colon #2 \rightarrow #3$}}
\newcommand{\dl}{{\cal D}_*}
\newcommand{\du}{{\cal D}^*}
\newcommand{\td}{{D}}

\begin{document}
\bibliographystyle{plain}
\title[The Combinatorics of Steenrod Operations]{The Combinatorics of Steenrod Operations on the Cohomology of Grassmannians }
\author{Cristian Lenart}
\address{Department of Mathematics, Massachusetts Institute of Technology, Cambridge, MA 02139}
\email{lenart@math.mit.edu}

\thanks{The author's visit to MSRI was supported by NSF grant DMS-9022140.}

\begin{abstract}
The study of the action of the Steenrod algebra on the mod $p$ cohomology of spaces has many applications to the topological structure of those spaces. In this paper we present combinatorial formulas for the action of Steenrod operations on the cohomology of Grassmannians, both in the Borel and the Schubert picture. We consider integral lifts of Steenrod operations, which lie in a certain Hopf algebra of differential operators. The latter has been considered recently  as a
realization of the Landweber-Novikov algebra in complex
cobordism theory; it also has connections with the action of the Virasoro algebra on
the boson Fock space. Our formulas for Steenrod operations are based on combinatorial methods which have not been used before in this area, namely  
Hammond operators and the combinatorics of Schur functions. 
We also discuss several applications of our formulas to the geometry of Grassmannians.
\end{abstract}

\maketitle

\section{Introduction}

The action of the {\em Steenrod algebra} on the mod $p$ cohomology of a space has a deep geometrical significance, whence it offers more information about the geometry of that space than just the cup product structure of its cohomology. Indeed, there are many instances where viewing the cohomology as a module over the Steenrod algebra distinguishes between spaces with isomorphic cohomology rings. Furthermore, the obstruction for maps induced in cohomology to be maps of algebras over the Steenrod algebra is a reasonably strong obstruction for the existence of certain topological maps. The Steenrod algebra can also be used to investigate the attaching maps of CW-complexes, as we discuss in \S3 and \S6. All these facts show the importance of determining the action of Steenrod operations on the cohomology of spaces we want to study. Among such spaces, the {\em Grassmannians} play an important role. Apart from projective spaces, very little is known about the attaching maps of their cells. 

From a combinatorial point of view, Grassmannians are important because their cohomology ring can be realized naturally in terms of symmetric functions, the Schur functions corresponding to the {\em Schubert classes}. No good explanation has been found yet for the occurence of Schur functions in both the cohomology of Grassmannians and the representation theory of the symmetric and general linear groups. However, there has been considerable interest recently in connections between the Steenrod algebra and the modular representation theory of the symmetric and general linear groups (see \cite{schums}, \cite{kuhgrf}, \cite{kstrsa}, \cite{cawpso}). From this perspective, it is possible that the Steenrod algebra might bring more evidence that the occurence of Schur functions in the two areas mentioned above is not a simple accident, and that there are deep connections between these areas. 

In this paper we present formulas for the action of Steenrod operations on the cohomology of Grassmannians, both in the Borel and the Schubert picture. These formulas have applications to the study of the attaching maps of the Schubert cells in the Grassmannians; they might also be relevant in the representation theoretical context discussed above. Our work is based on combinatorial methods which have not been used before in this area, namely  
Hammond operators and the combinatorics of Schur functions. 

In \S2, we present a Hopf algebra of differential operators, which has been considered recently by Reg Wood \cite{woodos} and Nigel Ray \cite{raypc} as a realization of the Landweber-Novikov algebra in complex cobordism, and which provides integral lifts of Steenrod operations. We point out that the Hopf algebra of symmetric functions has a natural {\em module Hopf algebra} structure over the algebra of differential operators. Some connections with the action of the {\em Virasoro algebra} on the {\em boson Fock} space are also mentioned. In \S3, we review the topological background concerning the Steenrod algebra and the cohomology of Grassmannians. In \S4, we present an efficient algorithm for computing the action of Steenrod operations on {\em Chern classes}. This leads to nice generalizations for primes $p>2$ of the Wu formula (for $p=2$); previous attempts to generalize it produced only partial results or very complicated closed formulas. In \S5, we investigate the action of the integral lifts of Steenrod operations introduced in \S2 on Schur functions; this corresponds via reduction mod $p$ to the action of Steenrod operations on Schubert classes in the cohomology of Grassmannians. We give several formulas involving the combinatorics of Young diagrams, analyze their consequences, and conjecture a generalization of one of them. In \S6 we apply our formulas in \S5 to some simple examples involving Schubert classes in Grassmannians, and suggest possible implications for the corresponding attaching maps. 

I am grateful to Nigel Ray and Reg Wood for introducing me to this area, and for pointing out to me the rich combinatorics involved. I am also grateful to Sergey Fomin, Mike Hopkins, Haynes Miller, and Gian-Carlo Rota for valuable discussions.

\section{A Hopf Algebra of Differential Operators}

In this section we present a Hopf algebra of differential operators, which has been considered recently by R. Wood \cite{woodos} as a framework for integral lifts of Steenrod operations. N. Ray \cite{raypc} showed that the algebra of differential operators is isomorphic to the Landweber-Novikov algebra in complex cobordism. Here, we postpone any reference to topology to \S3. We refer the reader to \cite{monhat} for all information concerning Hopf algebras.

All rings and algebras we consider in this paper are assumed
graded by complex dimension, so that products commute without signs. 

We start by recalling some concepts and notation from \cite{woodos}. The {\em Weyl algebra} $\calw^*$ is the associative algebra with unit (over the rationals) which is generated by $x_i$, $\partial_j$ ($i,j$ non-negative integers) subject to the relations
\[[x_i,x_j]=0\,,\spa[\partial_i,\partial_j]=0\,,\spa[x_i,\partial_j]=\delta_{i,j}\,;\]
here the square brackets denote the Lie product, and $\delta_{ij}$ the Kronecker delta. For every partition $\lambda\vdash n$, we let $x_\lambda:=x_{\lambda_1}x_{\lambda_2}\ldots$, and similarly for $\partial_\lambda$. In this paper, we use the standard notation for partitions of positive integers, namely given $\lambda=(\lambda_1,\lambda_2,\ldots,\lambda_l)=(1^{m_1}2^{m_2}\ldots k^{m_k})$ with $\lambda_1\ge\lambda_2\ge\ldots\ge\lambda_l>0$, we write $\lambda'$ for its {\em conjugate}, and
\[l(\lambda):=l\,,\quad |\lambda|:=\lambda_1+\ldots+\lambda_l\,,\quad \|\lambda\|:=m_1!\ldots m_k!\,;\]
recall that $l(\lambda)$ is known as the length of $\lambda$, and $|\lambda|$ as the weight of $\lambda$. Consider the graded polynomial algebra $W^*:=\bZ[x_1,x_2,\ldots]$ (every $x_i$ has degree 1). There is a natural action of the Weyl algebra on $W^*$ ($\partial_i$ acts as partial derivative with respect to $x_i$), and a natural grading ($x_\lambda\,\partial_\mu$ has degree $l(\lambda)-l(\mu)$). The Weyl algebras $\calw_n$ are defined for each $n$ in a similar way, by restricting to the finite set of variables $x_1,\ldots,x_n$ and the corresponding partial derivatives. For every integer $k$, we consider the summands $\calw_n^k$ of degree $k$, the inverse limit $\widehat{\calw}^k$ with respect to the restriction maps, and $\widehat{\calw}^*:=\bigoplus_k \widehat{\calw}^k$. Clearly, the action of $\calw^*$ on $W^*$ extends to an action of $\widehat{\calw}^*$. 

We now consider certain elements in $\widehat{\calw}^*$. Let $D_0$ be the unit of $\widehat{\calw}^*$, and let 
\[\td_k:=\sum_{i\ge 1} x_i^{k+1}\partial_i\]
for every $k\ge 1$. A new multiplication, denoted by $\vee$, is defined on $\widehat{\calw}^*$ in \cite{woodos} by
\[x_\lambda\partial_\mu\vee x_\nu\partial_\pi:=x_\lambda x_\nu\partial_\mu\partial_\pi\,;\]
this is a good definition, because every element of $\widehat{\calw}^*$ can be written uniquely in a form (known as the {\em standard form}) where the polynomial part is on the left, and the derivatives are on the right. For every partition $\lambda=(1^{m_1}2^{m_2}\ldots k^{m_k})$, we define
\[\td_\lambda:=\frac{\td_1^{\vee m_1}}{m_1!}\vee\ldots\vee \frac{\td_k^{\vee m_k}}{m_k!}\,,\]
where $\td_i^{\vee n}$ denotes the $n$-fold $\vee$-multiplication; here we regard $D_n$ as an abbreviation for $D_{(n)}$, and $D_0$ as the element $D_\lambda$ corresponding to the empty partition. We denote by $\du$ the $\bZ$-module spanned by all $\td_\lambda$. The following result is proved in \cite{woodos}.

\begin{thm}\cite{woodos}\label{lndiff}\hfill
\begin{enumerate}
\item $\du$ is closed under composition in $\widehat{\calw}^*$. Furthermore, $\du$ is a graded cocommutative Hopf algebra with the following comultiplication and counit:
\[\varDelta(D_\lambda)=\sum_{\mu\cup\nu=\lambda}D_\mu\otimes D_\nu\,,\qquad \casesc{\varepsilon(D_\lambda)}{=}{1}{\lambda=\emptyset}{0}\]
here the union of partitions is defined by concatenating them and reordering the parts. \item $W^*$ is a left {\em $\du$-module algebra}, i.e. the multiplication and the unit maps of $W^*$ are $\du$-module maps. In particular, the so-called Cartan formula holds
\begin{equation}\label{cartan} 
D_\lambda\,(yz)=\sum_{\mu\cup\nu=\lambda} (D_{\mu}\, y)\,(D_{\nu}\, z)\,,
\end{equation}
for all $y,z$ in $W_*$.
\end{enumerate}
\end{thm}

\begin{rem}\hfill{\rm
\begin{enumerate}
\item It is known from \cite{woodos} that the algebra $\du\otimes\bQ$  is generated under composition by $D_1$ and $D_2$ only.
\item It is easy to check that the Lie bracket of $D_k$ and $D_l$ is given by $[D_k,D_l]=(l-k)D_{k+l}$. Therefore, the $\bZ$-module spanned by $D_k$ is a Lie algebra, and its rational universal enveloping algebra is $\du\otimes\bQ$.
\item $D_\lambda$ has degree $|\lambda|$, and maps $W^n$ (that is polynomials of degree $n$) to $W^{n+|\lambda|}$. Hence the action map $\du\otimes W^*\rightarrow W^*$ is a map of graded algebras.
\item We have that $D_\lambda\,x_\mu=0$ whenever $l(\lambda)>l(\mu)$.
\item The action of $\du$ on $W^*$ is clearly faithful. On the other hand, the space of invariants of this action consists of $0$ only; this can be easily seen by considering the action of $D_k$ on a homogeneous polynomial $y$, and observing that the leading monomial of $D_k\,y$ can be expressed in terms of the leading monomial of $y$.
\end{enumerate}
}
\end{rem}

The non-commutative multiplication (composition) in $\du$ will be denoted by $\circ$. R. Wood \cite{woodos} derives a formula for expanding any composition $D_\lambda\circ D_\mu$ in the basis of $\du$. By duality, he also obtains the formula for the comultiplication in the graded dual $\dl$ of $\du$. Since the Hopf algebra structure of $\dl$ is interesting from a combinatorial point of view, we describe it here too. 

As an algebra, $\dl$ is just a graded polynomial algebra $\bZ[b_1,b_2,\ldots]$, where the element $b_i$ has degree $i$ for every $i$. The element $D_\lambda$ of $\du$ is dual to $b_\lambda:=b_{\lambda_1}b_{\lambda_2}\ldots$ with respect to the monomial basis of $\dl$. The non-cocommutative comultiplication in $\dl$ is given by the following formula: 
\begin{equation}\label{coprod}
\varDelta(b_k)=\sum_{i=0}^k (b)_{k-i}^{i+1}\otimes b_i\,,
\end{equation}
where $b_0:=1$, $b$ is the formal sum $b_0+b_1+b_2+\ldots$,  and  $(b)_{k-i}^{i+1}$ denotes the part of degree $k-i$ in $b^{i+1}$. In fact, this amounts to substitution of power series; indeed, if we let $b(t):=t+b_1 t^2+b_2 t^3+\ldots$, and we let $b^{(1)}(t)$, $b^{(2)}(t)$ be the corresponding generating functions for $b_k\otimes 1$ and $1\otimes b_k$, we have 
\begin{equation}\label{comulttwo}
\sum_{k\ge 0}\varDelta(b_k)t^{k+1}=b^{(2)}(b^{(1)}(t))\spa\spa\mbox{in}\spa(\dl\otimes\dl)[[t]]\,.
\end{equation}
The counit is specified by 
\[\casec{\varepsilon(b_k)}{=}{1}{k=0}{0}\]
for all $k\ge 0$. As discussed in \cite{schiha}, $\dl$ is the {\em incidence Hopf algebra} corresponding to the family of compositions of integers, which are Boolean algebras. Another combinatorial aspect concerning the structure of $\dl$ is the fact that its antipode is given by Lagrange inversion. More precisely, $S(b_k)$ is the coefficient of $t^{k+1}$ in the substitutional inverse of $b(t)$, that is the coefficient of $t^{k}$ in
\[\frac{1}{k+1}\left(\frac{t}{b(t)}\right)^{k+1}\,.\]
If we think of $b_k$ as being the complete homogeneous symmetric functions, then the antipode is precisely the second involution on symmetric functions considered in Macdonald's book (\cite{macsfh}, p. 35).

As pointed out in \cite{woodos}, the $\vee$-multiplication in $\du$ is the transpose of the comultiplication in $\dl$ specified by $b_k\mapsto1\otimes b_k+b_k\otimes 1$. This shows that $\du$ is also a Hopf algebra with respect to the $\vee$-multiplication and the same comultiplication. In fact, it is a tensor product of divided power Hopf algebras in $D_k$, for $k\ge 1$.

Let us now consider the algebra of symmetric functions with integer coefficients in $x_1,x_2,\ldots$, which we denote by $Sym_*$. We use the notation of \cite{macsfh} for symmetric functions, namely $h_\lambda$ for the complete homogeneous symmetric functions, $e_\lambda$ for the elementary symmetric functions, $p_\lambda$ for the power sums, $m_\lambda$ for the monomial symmetric functions, and $s_\lambda$ for the Schur functions. It is well-known that $Sym_*$ is a polynomial algebra (over the integers) in both $h_n$ and $e_n$, and that $Sym_*\otimes\bQ$ is a polynomial algebra (over the rationals) in $p_n$. In fact, $Sym_*$ is a Hopf algebra in which $p_n$ are primitives, and
\begin{equation}\label{comsym}
\delta(h_n)=\sum_{k=0}^n h_k\otimes h_{n-k}\,,\quad \delta(e_n)=\sum_{k=0}^n e_k\otimes e_{n-k}\,.
\end{equation}
The graded dual $Sym^*$ is isomorphic to $Sym_*$, and we have the standard pairing \linebreak $\dual{m_\lambda}{h_\mu}=\delta_{\lambda\mu}$. The transpose of the multiplication by a symmetric function $f$ is denoted by $f^\perp$, and the corresponding operator is known as a {\em Hammond operator}. For instance $p_n^\perp=n\,\partial/\partial p_n$. As pointed out in \cite{macsfh}, the linear span of $p_n$ and $p_n^\perp$ is a Heisenberg Lie algebra.

The action of $\du$ on $W^*$ can be extended to an action on the power series ring $\widehat{W}_*:=\bZ[[x_1,x_2,\ldots]]$, and then restricted to a left action on $Sym^*$. By partial duality, we obtain a right action of $\du$ on $Sym_*$, which we denote by $D^\perp\,f$ for every $f$ in $Sym_*$ and $D$ in $\du$. On the other hand, if we consider duality rather than partial duality, we obtain a left coaction of $\dl$ on $Sym_*$. In fact, there are two more actions and three more coactions that can be obtained in a similar way, and any of them determines all the others (see \cite{boaewb}).

According to Theorem \ref{lndiff} (2), $Sym^*$ is a left $\du$-module algebra. We now show that in fact it is a {\em module Hopf algebra}, which means that the comultiplication and the antipode are also $\du$-module maps. 

\begin{thm}
$Sym^*$ is a left $\du$-module Hopf algebra, and $Sym_*$ is a right $\du$-module Hopf algebra. On the other hand, $Sym_*$ is a left $\dl$-comodule Hopf algebra.
\end{thm}
\begin{proof}
We need to check that the so-called co-Cartan formula holds:
\begin{equation}\label{cocartan}
\delta(D\,f)=(\varDelta(D))\,(\delta(f))\,,
\end{equation}
where $D$ is an arbitrary element in $\du$, $f\in Sym^*$, $\varDelta$ is the comultiplication in $\du$, and $\delta$ is the comultiplication in $Sym^*$. This can be easily checked for $f=p_n$ and $D=D_\lambda$ using the fact that $D_\lambda\,p_n=\beta\, p_{n+|\lambda|}$ for a certain integer $\beta$ (cf. the computation in the proof of Proposition \ref{coactsym}). On the other hand, it is also easy to show that if (\ref{cocartan}) holds for the symmetric functions $f_1$ and $f_2$, it holds for $f_1f_2$, by using the Cartan formula (\ref{cartan}). A similar argument applies to the antipode. The other statements can be obtained by partial duality and duality, respectively.  
\end{proof} 

\newpage

\begin{rem}\hfill{\rm
\begin{enumerate}
\item The condition for the comultiplication in $Sym^*$ to be a $\du$-module map is equivalent to the map $\du\otimes Sym^*\rightarrow Sym^*$ being a coalgebra map. 
\item Similarly, the condition for the multiplication in $Sym_*$ to be a $\dl$-comodule map is equivalent to the coaction map $Sym_*\rightarrow \dl\otimes Sym_*$ being an algebra map.
\end{enumerate}
}
\end{rem}

The coaction of $\dl$ on $Sym_*$ is also related to substitution of power series. In order to express it, we introduce the generating function for the complete homogeneous symmetric functions $h(t):=1+h_1 t+h_2 t^2+\ldots$; we also let $h^{(2)}(t)$ be the corresponding generating function for $1\otimes h_k$  (note the difference from the definition of $b(t)$). 

\begin{prop}\label{coactsym}
The coaction map $Sym_*\rightarrow \dl\otimes Sym_*$ (which we observed that is an algebra map) is specified by
\[h(t)\mapsto h^{(2)}(b^{(1)}(t))\,.\]
\end{prop}
\begin{proof}
Let us denote by $\widetilde{\varDelta}$ the coaction map. The coefficient of $b_\lambda\otimes h_\mu$ in $\widetilde{\varDelta}(h_n)$ is
\[\alpha_{\lambda\mu}=\dual{D_\lambda\otimes m_\mu}{\widetilde{\varDelta}(h_n)}=\dual{D_\lambda\,m_\mu}{h_n}\,,\]
where $|\lambda|+|\mu|=n$. By the definition of $D_\lambda$, this coefficient is 0 unless $l(\mu)=1$. So let $\mu=(k)$, where $k:=n-|\lambda|$. By (\ref{cartan}), $D_\lambda \,x_i^k=\beta\,x_i^{k+|\lambda|}$, where $\beta$ is the number of terms of the form $D_{i_1}\otimes\ldots\otimes D_{i_{k}}$ in the $(k-1)$-fold comultiplication $\varDelta_{k-1}(D_\lambda)$. Assuming that $l(\lambda)\le k$, we can see that $\beta$ is the number of permutations of the multiset consisting of the parts of $\lambda$ and $k-l(\lambda)$ zeros. Hence
\[\case{\alpha_{\lambda(k)}=\beta}{=}{\frac{k!}{\|\lambda\|(k-l(\lambda))!}}{l(\lambda)\le k}{0}\]
On the other hand, by the multinomial theorem, this is exactly the coefficient of $b_\lambda\, t^{|\lambda|+k}$ in $b(t)^k$. This concludes the proof.
\end{proof}

Finally, we show that the actions of the operators $D_k$ and $D_k^\perp$ on $Sym^*$ and $Sym_*$ are closely related to the action of the {\em Virasoro algebra} on the {\em boson Fock space}.

\begin{prop}\label{expd}
For every $k\ge 1$, the operators $D_k$ and $D_k^\perp$ on $Sym^*$ and $Sym_*$ can be expressed as
\[ D_k=\sum_{i\ge 1}     p_{i+k}\,p_i^\perp\,,\quad D_k^\perp=\sum_{i\ge 1}     p_i\,p_{i+k}^\perp\,.\]
\end{prop}
\begin{proof}
The actions of $D_k$ and $\sum_{i\ge 1}     p_{i+k}\,p_i^\perp$ clearly agree on $p_j$, since $D_k\,p_j=jp_{j+k}$. On the other hand, both operators are derivations, whence their actions agree on any symmetric function. The proof of the statement concerning $D_k^\perp$ is similar.
\end{proof}

We briefly recall some basic concepts from conformal field theory in order to discuss the connection with the algebra $\du$; for more details, the reader is referred to \cite{karhwr}, for instance. The boson Fock space of central charge $1$ is defined to be
\[\calh_T:=\bQ[p_1,p_2,\ldots][u,u^{-1}]\,,\]
where $p_k$ can be thought of as the power sum symmetric functions. The charge $n$ sector of $\calh_T$ is $\bQ[p_1,p_2,\ldots]u^n$. The Virasoro algebra of central charge $c$ is the Lie algebra over $\bQ$ spanned by elements $L_k$, $k\in \bZ$, which has the following bracket:
\[[L_k,L_l]=(k-l)L_{k+l}+\frac{k^3-k}{12}\,c\,\delta_{k,-l}\,.\]
There is a natural action of the Virasoro algebra on the boson Fock space. For $c=1$, the basis elements of the Virasoro algebra act on the charge 0 sector of $\calh_T$ in the following way:
\begin{gather*}
L_0=\sum_{i\ge 1}     p_i\,p_i^\perp\,,\\
L_k=\sum_{i\ge 1}     p_i\,p_{i+k}^\perp+\frac{1}{2}\sum_{i=1}^{k-1}p_i^\perp\,p_{k-i}^\perp\,,\\
L_{-k}=\sum_{i\ge 1}     p_{i+k}\,p_{i}^\perp+\frac{1}{2}\sum_{i=1}^{k-1}p_i\,p_{k-i}\,,
\end{gather*}
where $k\ge 1$. We immediately observe that the expressions for $D_k$ and $D_k^{\perp}$ in Proposition \ref{expd} are just truncated versions of the expressions for $L_{-k}$ and $L_k$, respectively (in fact $D_1=L_{-1}$ and $D_1^\perp=L_1$). Furthermore, the Lie algebras spanned by $\setp{D_k}{k\ge 0}$ and $\setp{D_k^\perp}{k\ge 0}$ are isomorphic to the corresponding Lie subalgebras of the Virasoro algebra. However, the elements $D_k$ and $D_k^\perp$ considered together do not span a Lie algebra.

We conclude by mentioning that T. Katsura et al. \cite{ksuccr} constructed a natural embedding of the dual of the Landweber-Novikov algebra (which is just the Hopf algebra $\dl$, as shown in the next section) in the charge 0 sector of the boson Fock space.

\section{Background from Algebraic Topology}

In this section, we briefly review the topology underlying the algebra presented in \S2. We refer the reader to \cite{saeco} and \cite{milsad} for all information concerning the Steenrod algebra, to \cite{adashg} for information about the Landweber-Novikov algebra, and to \cite{fulyt} or \cite{hilgcg} for information about the cohomology of Grassmannians.

For every prime $p$, the mod $p$ {\em Steenrod algebra} ${\cal A}(p)$ is the graded associative algebra of mod $p$ stable operations in ordinary cohomology theory over ${\mathbb F}_p$. If $p=2$, it is generated by the Steenrod squares $Sq^n$ ($n\ge 1$), modulo the Adem relations; if $p$ is odd, it is generated by the Bockstein operation and the (reduced) Steenrod operations ${\cal P}^n$ ($n\ge 1$), modulo the corresponding Adem relations. The operations $Sq^n$ raise degree by $n$, while the operations ${\cal P}^n$ raise degree by $2n(p-1)$; in other words:
\[\map{Sq^n}{H^k(-,{\mathbb F}_p)}{H^{k+n}(-,{\mathbb F}_p)}\,,\quad\map{{\cal P}^n}{H^k(-,{\mathbb F}_p)}{H^{k+2n(p-1)}(-,{\mathbb F}_p)} \,.\]
For simplicity, we will denote $Sq^{2n}$ by ${\cal P}^n$ as well. Milnor showed that ${\cal A}(p)$ modulo the two-sided ideal generated by the Bockstein operation has a linear basis indexed by compositions of positive integers. Furthermore, he showed that both ${\cal A}(p)$ and its dual have a natural Hopf algebra structure over ${\mathbb F}_p$. The action of the Steenrod algebra on the cohomology of spaces gives rise (by partial duality) to an action of it on homology, and to a coaction of the dual Steenrod algebra on homology (by duality). There are two more actions and three more coactions that can be obtained in a similar way, and they are discussed in \cite{boaewb}.

There are many generalized cohomology theories now used in algebraic topology, and most of them have a Hopf algebra of stable operations analogous to the Steenrod algebra. In every such situation, there is a geometric way to define four actions and four coactions of these algebras or their duals on the corresponding cohomology or homology. An important such cohomology theory is complex cobordism $MU^*(-)$ introduced by Milnor in \cite{milcro}. The structure of the algebra $MU^*(MU)$ of all operations in complex cobordism was determined by Landweber \cite{lancoh} and Novikov \cite{novmat}. They showed that $MU^*(MU)$ is isomorphic to the tensor product of the complex cobordism ring $MU_*$ and the so-called {\em Landweber-Novikov algebra}, which turns out to be precisely the Hopf algebra $\du$ discussed in the previous section; in fact $MU^*(MU)$ is a so-called {\em Hopf algebroid}.

Let us now refer to the Grassmannian of $n$-planes in $\bC^\infty$, which is also the classifying space for principle $U(n)$-bundles, whence the notation $BU(n)$; note that $BU(1)$ is the infinite complex projective space, usually denoted by $\CP$. There is a natural inclusion of $BU(n)$ into $BU(n+1)$, so we can consider the union of all $BU(n)$, which is the space denoted by $BU$. Let $E^*(-)$ be an unreduced multiplicative cohomology theory with complex orientation $Z\in E^2(\bC P^\infty)$; in particular, $E^*(-)$ can be ordinary cohomology (with integer or ${\mathbb F}_p$ coefficients) or complex cobordism. We write $E_*:=E_*(\mbox{point})$ for the coefficient ring, and identify $E_n$ with $E^{-n}:=E^{-n}(\mbox{point})$, as it is usually done. It is well-known that $E^*(BU(n))\cong E^*[c_1,c_2,\ldots, c_n]$ and $E^*(BU)\cong E^*[c_1,c_2,\ldots]$, where $c_i$ are the generalized {\em Chern classes}. This description of the cohomology of $BU(n)$ is called the {\em Borel picture}. It is also known that the map
\[E^*(BU(n))\rightarrow E^*(\CP\times\ldots\times\CP)\cong E^*[x_1,x_2,\ldots,x_n]\]
induced by the classifying map of the direct product of $n$ copies of
the Hopf bundle over $\bC P^\infty$ is a monomorphism mapping $c_i$, with $i\le n$, to
the $i$-th elementary symmetric polynomial in $x_1:=Z\otimes 1\otimes\ldots\otimes
1$, $x_2:=1\otimes Z\otimes\ldots\otimes 1$, ..., $x_n:=1\otimes
1\otimes\ldots\otimes Z$. On the other hand, we have that
$E_*(BU)\cong E_*[b_1,b_2,\ldots]$, and that $c_n$ is dual to $b_1^n$
with respect to the monomial basis of $E_*[b_1,b_2,\ldots]$. The
multiplicative structure of $E_*(BU)$ is determined by the map
$BU\times BU\rightarrow BU$ classifying the Whitney sum of vector
bundles. The diagonal map $BU\rightarrow BU\times BU$ induces a comultiplication $\delta\colon E_*( BU)\rightarrow E_*( BU\times BU)\cong E_*( BU)\otimes E_*( BU)$ satisfying
\[\delta(b_n)=\sum_{i=0}^n b_i\otimes b_{n-i}\,,\]
which turns $E_*( BU)$ into a Hopf algebra. It follows from the above considerations that we may identify
$E_*(BU)$ with $Sym_*\otimes E_*$ and $E^*(BU)$ with $Sym^*\otimes E^*$, in such a way
that $b_n$ is identified with $h_n$ and $c_n$ with $e_n$. 

We can also consider the Grassmannian $Gr_n(\bC^{n+k})$ of $n$-planes in $\bC^{n+k}$, whose cohomology ring is a certain quotient of the cohomology of $BU(n)$; indeed, the obvious inclusion $Gr_n(\bC^{n+k})\hookrightarrow BU(n)$ induces the quotient map. Let $\lambda$ be an arbitrary partition (identified with its Young diagram) with at most $n$ rows and $k$ columns, and let $F_\bullet$ be a fixed complete flag
\[0=F_0\subset F_1\subset F_2\subset \ldots \subset F_{n+k}=\bC^{n+k}\]
of subspaces with $\mbox{dim}(F_i)=i$. There is a {\em Schubert cell} $\varOmega_\lambda=\varOmega_\lambda(F_\bullet)$ defined by
\[\varOmega_\lambda:=\setp{V\in Gr_n(\bC^{n+k})}{\mbox{dim}(V\cap F_{i+\lambda_{n+1-i}})\ge i\,,\spa1\le i\le n}\,;\]
here we set $\lambda_i=0$ for $i>l(\lambda)$. It is known that $\varOmega_\lambda$ is an irreducible closed subvariety of $Gr_n(\bC^{n+k})$ of dimension $2|\lambda|$. It is also known that $\varOmega_\mu\subset\varOmega_\lambda$ if and only if the Young diagram of $\mu$ is strictly contained in that of $\lambda$, and that $\varOmega_\lambda$ is the disjoint union of the interiors of the cells $\varOmega_\mu$, $\mu\subseteq\lambda$. Furthermore, the cohomology class $\sigma_\lambda$ of $\varOmega_\lambda$ (that is the dual of the fundamental class of $\varOmega_\lambda$) is independent of the choice of the fixed flag defining it. These classes, known as Schubert classes, form a basis over the integers for the cohomology ring of $Gr_n(\bC^{n+k})$. This description of the cohomology of $Gr_n(\bC^{n+k})$ is called the {\em Schubert picture}. The crucial fact is that the quotient map from the cohomology of $BU(n)$ (which we identified with the ring of symmetric polynomials in $n$ variables) to the cohomology ring of $Gr_n(\bC^{n+k})$ maps the Schur polynomial $s_\lambda(x_1,\ldots,x_n)$ to $\sigma_\lambda$ if the Young diagram $\lambda$ has at most $n$ rows and $k$ columns, and to 0 otherwise. 

Another topological fact we need concerns the relation between the action of Steenrod operations on Schubert cells and the corresponding attaching maps. It follows from the naturality of Steenrod operations that any composition of them applied to a Schubert class is a linear combination (with non-zero coefficients) of Schubert classes whose corresponding cells are non-trivially attached to the initial cell. Hence, whenever we detect a non-zero coefficient, we can say that we have a non-trivial attaching map. Conversely, when every composition of Steenrod operations applied to a Schubert class gives a zero coefficient for a higher dimensional Schubert class, we have a chance that the corresponding attaching map is trivial. 
 
Let us now consider the infinite product of infinite complex projective spaces, whose $E^*$-cohomology is $W^*\otimes E^*=E^*[x_1,x_2,\ldots]$. The action of the Landweber-Novikov algebra on $MU^*[x_1,x_2,\ldots]$ restricts to an action on $W^*$, which coincides with the action of the algebra $\du$ discussed in \S2. This was proved by N. Ray \cite{raypc} by comparing the action of the operators $D_\lambda$ with the action of the Landweber-Novikov operations described, for instance, in \cite{adashg}. On the other hand, the action of the Steenrod operations on $\bF_p[x_1,x_2,\ldots]$ is given by the following rules:
\begin{enumerate}
\item $\calp^0$ is the identity;
\item $\calp^1\,x_i=x_i^p$ and $\calp^n\,x_i=0$ for $n>1$ ($i\ge 1$);
\item $\calp^n\,(yz)=\sum_{k=0}^n (\calp^k\,y)\,(\calp^{n-k}\,z)$ for every $y,z$ in $W^*\otimes\bF_p$ (Cartan formula).
\end{enumerate}
In fact, the Cartan formula holds for every space, while $\calp^n$ raises every cohomology class $x$ to its $p$-th power if dim$(x)=2n$, and maps it to 0 if dim$(x)<2n$. This justifies the name ``reduced power operations'' (or ``squaring operations'' for $p=2$). As far as the action on the cohomology of $BU$ is concerned, we have the following formula due to Borel and Serre \cite{basglp}:
\begin{equation}\label{stchern}
\calp^n\,c_{k+n}=m_{(1^kp^n)}\,;
\end{equation}
here $m_\lambda$ is the corresponding monomial symmetric function in $H^*(BU,\bF_p)$, identified with the ring of symmetric functions with $\bF_p$ coefficients. This formula was used in \cite{lansdl} to lift the action of $\calp^n$ on the mod $p$ cohomology of $BU$ to the integral cohomology. It is not difficult to see that the operator defining this integral lift is precisely the differential operator $D_{((p-1)^n)}$ in $\du$. In fact, as the three rules above show, this operator lifts the action of $\calp^n$ on the mod $p$ cohomology of the infinite product $\CP\times\CP\times\ldots$ to the integral cohomology. Furthermore, the subalgebra of $\du$ generated by the operators $D_{((p-1)^n)}$ ($n\ge 0$) under composition can be identified upon reduction mod $p$ with the quotient of ${\cal A}(p)$ by the two-sided ideal generated by the Bockstein operation. This is proved in \cite{woodos} for $p=2$, but the argument can be easily extended to any prime. Let us mention R. Wood's notation $SQ^n$ for $D_{(1^n)}$, which we will also use.

\section{Steenrod Operations in the Borel Picture of the Cohomology of Grassmannians}

In this section we present an efficient algorithm for computing the action of Steenrod operations on Chern classes in $H^*(BU,\bF_p)$. This problem has been studied for a long time. According to (\ref{stchern}), it is just the problem of expanding the monomial symmetric function $m_{(1^kp^n)}$ in the basis of elementary symmetric functions. We can ask for this expansion over the integers, which would correspond to the integral lift of Steenrod operations, or we could just ask for the mod $p$ reduction. In principle, we can then compute the action of Steenrod operations on any monomial in the Chern classes using the Cartan formula. 

For $p=2$, the mod 2 formula was found by Wu \cite{wucvg}: 
\begin{equation}\label{wu}
m_{(1^k2^n)}=\sum_{i=0}^n \binom{n+k-i-1}{n-i}\,e_i\,e_{2n+k-i}\quad\mbox{(mod 2)}\,.
\end{equation}
The corresponding integral formula was given by Carlitz \cite{carsif}. As far as the case $p>2$ is concerned, it is shown in \cite{lansdl} that the coefficient of $e_{pn+k}$ in $m_{(1^kp^n)}$ is $\frac{pn+k}{n+k}\binom{n+k}{n}$. This agrees with the mod $p$ result of Brown and Peterson \cite{bapsrs}, who show that the coefficient is $\binom{n+k-1}{n}$ mod $p$. Other partial results appeared in \cite{bpdhbs} and \cite{petmpw}. On the other hand, Shay \cite{shampw} obtains a closed formula for the expansion of $m_{(1^kp^n)}$ in the elementary symmetric functions by using the expansion in the power sums as an intermediate step. However, his formula is complicated and difficult to work with; indeed, it is a non-trivial matter to deduce from this formula when a given coefficient is 0. 

The general problem of expanding a monomial symmetric function $m_\lambda$ in the basis of elementary ones goes back to Cayley and MacMahon \cite{macca}. Their idea is to apply the Hammond operators $h_i^\perp$ to the
expansion of $m_\lambda$ with indeterminate coefficients. It is easy to see that $h_k^\perp$ acts on a monomial symmetric function by removing a part equal to $k$, if it exists, and gives 0 otherwise; furthermore, $h_k^\perp$ satisfies the Cartan formula. Hence we obtain a system of equations which depends on the expansion of monomial symmetric functions indexed by partitions obtained from $\lambda$ by removing one part. In this section, we discuss a refinement of this algorithm which is specially designed for the case $\lambda=(1^kp^n)$, and is based on symmetric functions in two sets of variables. Let us also note that there are combinatorial formulas for the expansion of an arbitrary monomial symmetric function in terms of the lattice of partitions of a set \cite{doufct7}, as well as in terms of certain combinatorial objects called {\em layered primitive bi-brick permutations} \cite{karpbb}. However, both of these formulas are complicated and difficult to work with, even in the special case we consider.

For the rest of this section, we fix an arbitrary integer $p\ge 2$ (not necessarily a prime). We use the following notation:
\begin{equation}\label{not}
m_{(1^kp^n)}=\sum_{\lambda\vdash pn+k}\alpha_\lambda^{n,k}\,e_\lambda\,.
\end{equation}
We also denote by $\lambda-1$ the partition obtained from $\lambda$ by subtracting $1$ from every part. The following Proposition collects a couple of simple observations.

\begin{prop}\label{trans}\hfill
\begin{enumerate}
\item We have $\alpha_\lambda^{n,k}=0$ unless $\lambda\ge (n^{p-1}(n+k))$ in the dominance order; furthermore, $\alpha_{(n^{p-1}(n+k))}^{n,k}=1$. In particular, $\alpha_\lambda^{n,k}=0$ unless $l(\lambda)\le p$ and $\lambda_1\ge n+k$.
\item If $l(\lambda)=p$, then $\alpha_{\lambda}^{n,k}=\alpha_{\lambda-1}^{n-1,k}$.
\end{enumerate}
\end{prop}
\begin{proof}
1) This result is in fact more general. We show first that for any $\lambda$, we have
\[e_\lambda=m_{\lambda'}+\sum_{\mu<\lambda'}\gamma_{\lambda\mu}\,m_\mu\,.\]
It is easy to see that the partitions $\mu$ occuring in the above sum are obtained in the following way: we place $\lambda_1'$ boxes in a column, then add $\lambda_2'$ boxes such that no two are in the same row, then continue in the same way. Clearly, $\mu\le\lambda'$ in the dominance order. It is also easy to see that the coefficient of $m_{\lambda'}$ is 1. We now use the above relation and induction on $\lambda$ to show that 
\[m_\lambda=e_{\lambda'}+\sum_{\nu>\lambda'} \theta_{\lambda\nu}\,e_\nu\,.\]
Here we also need the well-known fact that $\mu<\lambda$ if and only if $\mu'>\lambda'$.

2) This follows easily by applying the Hammond operator $h_p^\perp$ to both sides of (\ref{not}). 
\end{proof}

According to the previous observations, it suffices to determine the coefficients $\alpha_\lambda^{n,k}$ corresponding to partitions $\lambda$ of length at most $p-1$. The crucial ingredient for our algorithm is working with symmetric functions in two sets of variables $x=(x_1,x_2,\ldots)$ and $y=(y_1,y_2,\ldots)$. Let
\[h_n:=h_n(x)\,,\quad \overline{h}_n:=h_n(y)\,,\quad H_n:=h_n(x,y)\,,\]
and similarly for the elementary symmetric functions. There are corresponding Hammond operators $H_n^\perp$, and they also satisfy the Cartan formula
\[H_n^\perp\,(uv)=\sum_{k=0}^n (H_k^\perp\,u)\,(H_{n-k}^\perp\,v)\,.\]
It is known that $H_n=\sum_{k=0}^n h_k\,\overline{h}_{n-k}$.
Let us also note that $H_1^\perp\,e_n=e_{n-1}$ and $H_1^\perp\,\overline{e}_n=\overline{e}_{n-1}$.

Now let us fix $n$. We define $\overline{m}_{(1^kp^n)}$ in $\bZ[x,y]$ as follows:
\[\casel{\overline{m}_{(1^kp^n)}}{:=}{\sum_{\lambda\vdash pn+k}\alpha_\lambda^{n,k}\,e_{\lambda\setminus\lambda_1}\,\overline{e}_{\lambda_1}}{k\ge n(p-2)}{(H_1^\perp)^{n(p-2)-k}\,\overline{m}_{(1^{n(p-2)}p^n)}}\]
the notation $\lambda\setminus i$ means the partition obtained from $\lambda$ by removing the part $i$. Let us now introduce the notation
\[\overline{m}_{(1^kp^n)}=\sum_{\lambda}\beta_\lambda^{n,k}\,e_\lambda\,\overline{e}_{pn+k-|\lambda|}\,;\]
clearly, $\overline{m}_{(1^kp^n)}$ has this form. 

\begin{lem}\label{lemalgor} 
We have $H_1^\perp\,\overline{m}_{(1^kp^n)}=\overline{m}_{(1^{k-1}p^n)}$ for $k\ge 1$, $H_p^\perp\,\overline{m}_{(1^kp^n)}=\overline{m}_{(1^{k}p^{n-1})}$ for $n\ge 1$, and $H_i^\perp\,\overline{m}_{(1^kp^n)}=0$ for $i\not\in\{1,p\}$.
\end{lem}
\begin{proof}
To prove the first formula, it is enough to consider the case $k>n(p-2)$. By using Proposition \ref{trans}, we deduce that in this case $\lambda_1>\lambda_2$ for every partition $\lambda$ with at least two parts for which $\alpha_\lambda^{n,k}\ne 0$. Indeed, on the one hand we have $n+k>np-n$, and on the other hand we have $\lambda_1\ge n+k$, which implies $\lambda_2\le np-n$. We conclude that the terms of $H_1^\perp\,\overline{m}_{(1^kp^n)}$ have the form $\beta\,e_{\lambda\setminus\lambda_1}\,\overline{e}_{\lambda_1}$. The first formula now follows by combining this remark with the definition of $\overline{m}_{(1^{k-1}p^n)}$ and the fact that $H_1^\perp$ lifts the action of $h_1^\perp$ on $Sym_*$ to the ring $\bZ[e_1,e_2,\ldots;\overline{e}_1,\overline{e}_2,\ldots]$. The other two formulas follow in a similar way if $k>n(p-2)$. Otherwise, we use the commutativity of the operators $H_i^\perp$:
\begin{align*}
H_p^\perp\,\overline{m}_{(1^kp^n)}&=H_p^\perp((H_1^\perp)^{n(p-2)-k+1}\,\overline{m}_{(1^{n(p-2)+1}p^n)})=(H_1^\perp)^{n(p-2)-k+1}(H_p^\perp\,\overline{m}_{(1^{n(p-2)+1}p^n)})\\
&=(H_1^\perp)^{n(p-2)-k+1}\,\overline{m}_{(1^{n(p-2)+1}p^{n-1})}=\overline{m}_{(1^{k}p^{n-1})}\,;\end{align*}
a similar argument applies to $H_i^\perp$ with $i\not\in\{1,p\}$.
\end{proof}

The following Theorem contains the results underlying our algorithm.

\begin{thm}\label{algor}\hfill
\begin{enumerate}
\item The following recurrence relation holds for all $\lambda$ and $n,k\ge 1$:
\[\beta_\lambda^{n,k}-\beta_\lambda^{n,k-1}=(-1)^{p-1}\beta_\lambda^{n-1,k}\,.\]
\item We have $\beta_\lambda^{n,k}=0$ whenever $l(\lambda)>p-1$ or $|\lambda|>n(p-1)$. Furthermore, $\beta^{n,0}=0$ (mod $p$) if $n>\lambda_1$. 
\item We have $\alpha_\lambda^{n,k}=\sum_i\beta_{\lambda\setminus i}^{n,k}$, where the summation ranges over all distinct parts of $\lambda$. If $n+k>\lambda_1$, then $\alpha_{\lambda\cup(pn+k-|\lambda|)}^{n,k}=\beta_\lambda^{n,k}$.
\end{enumerate}
\end{thm}
\begin{proof} 1) We apply the operators $H_i^\perp$, $i=1,2,\ldots,p$ to $\overline{m}_{(1^{k}p^{n})}$, and identify the coefficient of $e_\lambda\,\overline{e}_{pn+k-i}$ in the result. By Lemma \ref{lemalgor}, we have
\begin{align*}
\beta_\lambda^{n,k-1}&=\beta_\lambda^{n,k}+\sum_{\lambda^{(1)}}c_{\lambda^{(1)}}\,\beta_{\lambda^{(1)}}^{n,k}\,,\\
0&=\sum_{\lambda^{(1)}}c_{\lambda^{(1)}}\,\beta_{\lambda^{(1)}}^{n,k}+\sum_{\lambda^{(2)}}c_{\lambda^{(2)}}\,\beta_{\lambda^{(2)}}^{n,k}\,,\\
&\ldots\\
0&=\sum_{\lambda^{(p-2)}}c_{\lambda^{(p-2)}}\,\beta_{\lambda^{(p-2)}}^{n,k}+\beta_{\lambda^{(p-1)}}^{n,k}\,,\\
\beta_\lambda^{n-1,k}&=\beta_{\lambda^{(p-1)}}^{n,k}\,.
\end{align*}
Here $\lambda^{(i)}$ ranges over partitions with at most $p-1$ parts obtained by adding $i$ boxes to the Young diagram of $\lambda$, such that no two of them are on the same row; $c_{\lambda^{(i)}}$ are certain positive integer coefficients which we do not specify. The recurrence relation for $\beta_\lambda^{n,k}$ now follows by taking the alternating sum of the formulas above.

2) Clearly, if $\beta_\lambda^{n,k}\ne 0$, then $l(\lambda)\le p-1$ by Proposition \ref{trans}. If $k\ge n(p-2)$, then $|\lambda|>n(p-1)$ implies $\beta_\lambda^{n,k}=0$, by the same Proposition; the definition of $\overline{m}_{(1^{k}p^{n})}$ now shows that this is true in general. If $n>\lambda_1$, then $\beta_\lambda^{n,0}=\alpha_{\lambda\cup(pn-|\lambda|)}^{n,0}$ (by assertion 3). But $\alpha_{\lambda\cup(pn-|\lambda|)}^{n,0}=0$ (mod $p$), since $e_n^p= m_{(p^n)}$ (mod $p$). 

3) This formula follows immediately from the fact that $H_1^\perp$ lifts the action of $h_1^\perp$ on $Sym_*$ to the ring $\bZ[e_1,e_2,\ldots;\overline{e}_1,\overline{e}_2,\ldots]$. If $n+k>\lambda_1$, then $\beta_{(\lambda\setminus i)\cup(pn+k-|\lambda|)}^{n,k}=0$ for all parts $i$ of $\lambda$ (by assertion 2). Hence, in this case $\alpha_{\lambda\cup(pn+k-|\lambda|)}^{n,k}=\beta_\lambda^{n,k}$. 
\end{proof}


The reason we are working in $\bZ[e_1,e_2,\ldots;\overline{e}_1,\overline{e}_2,\ldots]$ is that the recurrence relation in Theorem \ref{algor} (1) does not hold for the coefficients $\alpha_\lambda^{n,k}$. The main advantage of our method is that we can express the coefficient $\alpha_\lambda^{n,k}$ mod $p$ in a very nice form, as the following Corollary shows. Here, as well as throughout this paper, we adopt the usual conventions $\binom{i}{j}=0$ if $j>i$ or $j<0$. 

\begin{cor}\label{simpleex}
Let $p$ be a prime, $\lambda$ a partition with at most $p-1$ parts, and $n_0:=\lceil |\lambda|/(p-1)\rceil,\,n_1:=\lambda_1$. For all $n,k\ge 0$, we can express $\beta_{\lambda}^{n,k}$ mod $p$ in the following way:
\[\sum_{i=n_0}^{n_1} \beta_\lambda^{i,0}\binom{n+k-i-1}{n-i}\,.\]
\end{cor}
\begin{proof}
By the recurrence relation in Theorem \ref{algor} (1), we have 
\[(-1)^{p-1}\beta_\lambda^{n,k}=\beta_\lambda^{n,0}+\sum_{j=1}^k\beta_\lambda^{n-1,j}\,.\]
The expression for $\beta_\lambda^{n,k}$ immediately follows by induction using the fact that $\beta_\lambda^{n,0}=0$ (mod $p$) unless $n_0\le n\le n_1$ (cf. Theorem \ref{algor} (2)), and the iterated Pascal formula
\[\sum_{j=1}^k\binom{q+j-1}{q}=\binom{q+k}{q+1}\,.\]
\end{proof}

Note the similarity between our expression for $\beta_\lambda^{n,k}$ mod $p$ and the Wu formula. For $p=2$ we have $\beta_i^{n,0}=0$ if $n<i$, and $\beta_i^{n,0}=1$; if $n>i$ we have $\beta_i^{n,0}=0$ (mod 2). The Corollary gives $\beta_i^{n,k}=\binom{n+k-i-1}{n-i}$. But $\alpha_{(2n+k-i,i)}^{n,k}=\beta_i^{n,k}$ if we assume $2n+k-i\ge i$; indeed, if strict inequality holds, we have $2n+k-i>n$, so $\beta_{2n+k-i}^{n,k}=0$. Hence we obtain the Wu formula (\ref{wu}). Thus, our Corollary turns out to be the natural generalization of this formula to primes $p>2$. 

We now give recurrence relations for the coefficients $\beta_\lambda^{n,0}$, so that we have a complete description of an algorithm for expanding $m_{(1^kp^n)}$ in the basis of elementary symmetric functions. If we are interested in the expansion over the integers, we need all $\beta_\lambda^{n,0}$ for $n\ge n_0$, but if we are working mod $p$ (as is the case with Steenrod operations), we only need $\beta_\lambda^{n,0}$ for $n_0\le n\le n_1$. We then determine $\beta_\lambda^{n,k}$ for $k>0$, and $\alpha_\lambda^{n,k}$ as a certain sum of the former (cf. Theorem \ref{algor} (3)). Let us introduce the following notation: given two partitions $\lambda=(1^{l_1}\ldots k^{l_k})$ and $\mu=(1^{m_1}\ldots k^{m_k})$, we let
\[c_{\lambda\mu}:=\binom{m_{k-1}}{l_k-m_k}\binom{m_{k-2}}{l_{k-1}+l_k-m_{k-1}-m_k}\ldots\binom{m_1}{l_2+\ldots+l_k-m_2-\ldots-m_k}\,.\]

\begin{prop}\label{recur}
 For every partition $\lambda$, we have $\beta_\lambda^{n,0}=\beta_{\lambda-1}^{n-1,0}$ if $l(\lambda)=p-1$, and
\[\sum_\mu c_{\lambda\mu}\,\beta_\mu^{n,0}+\sum_\nu c_{\lambda\nu}\,\beta_\nu^{n,0}=0\,;\]
here $\mu$ ranges over partitions obtained from $\lambda$ by decrementing at most $p-1-l(\lambda)$ parts different from $1$ by $1$ and adding as many parts equal to 1; on the other hand, $\nu$ ranges over partitions obtained from $\lambda$ by decrementing $k\le p-2-l(\lambda)$ parts different from $1$ by $1$ and adding $k+1$ parts equal to $1$.
\end{prop}
\begin{proof}
The recurrence relation follows easily by applying the operator $H_{l(\lambda)+1}^\perp$ to $\overline{m}_{(1^kp^n)}$ and identifying the coefficient of $e_{\lambda-1}\,\overline{e}_{pn+k-|\lambda|-1}$ in the result. The coefficient $c_{\lambda\mu}$ represents the number of ways to obtain the partition $\lambda-1$ from $\mu$ by subtracting $1$ from certain parts of $\mu$, and then reordering the parts (see the diagram below, where $\lambda=(1^{l_1}\ldots k^{l_k})$ and $\mu=(1^{m_1}\ldots k^{m_k})$; the black boxes are removed from $\lambda$, and the bold ones are added in order to obtain $\mu$).

\[
\begin{array}{c}
\mbox{\psfig{file=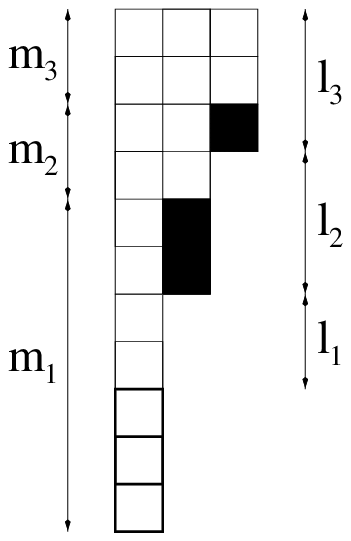}}
\end{array}
\]
\end{proof}

Note that $c_{\lambda\lambda}=1$, and that the above recurrence relations express $\beta_\lambda^{n,0}$ with $l(\lambda)<p-1$ in terms of $\beta_\mu^{n,0}$ with $l(\mu)>l(\lambda)$.

\begin{exa}\label{exaone}{\rm
We can use the above results to give a closed mod $p$ formula for $\alpha_\lambda^{n,k}$ in the case $p=3$. A similar formula was obtained in \cite{sugwfm}, by using the method in \cite{shampw}.

It follows from the discussion above that it suffices to consider partitions $\lambda$ with two parts. In other words, it suffices to compute $\beta_i^{n,k}$, which in turn depend on $\beta_i^{n,0}$. We know that $\beta_i^{n,0}=0$ if $i>2n$, and that $\beta_i^{n,0}=0$ (mod 3) if $n>i$. By Proposition \ref{recur}, we have $\beta_i^{n,0}=-\beta_{i-1}^{n-1,0}-\beta_{i-2}^{n-1,0}$ (we let $\beta_j^{n,0}=0$ if $j\le 0$). Hence $\beta_i^{n,0}$ mod 3 is completely determined by $\beta_1^{1,0}=-1$ and $\beta_2^{1,0}=-2$ (see the discussion in the next paragraph). A straightforward iteration of the recurrence relation gives
\[\beta_i^{n,0}=(-1)^j\sum_{l=0}^j\binom{j}{l}\beta_{i-j-l}^{n-j,0}\,.\]
Letting $j=n-1$, we obtain $\beta_i^{n,0}=(-1)^{n-1}\frac{2i}{n}\binom{n}{i-n}$ (mod 3). Combining this formula with Corollary \ref{simpleex} gives
\begin{equation}
\beta_i^{n,k}=\sum_{j=\lceil i/2\rceil}^i (-1)^{j-1}\frac{2i}{j}\,\binom{j}{i-j}\,\binom{n+k-j-1}{n-j}\quad \mbox{(mod 3)}\,.
\end{equation}
Finally, assuming that $3n+k\ge 2i$, we have 
\begin{equation}
\case{\alpha_{(3n+k-i,i)}^{n,k}}{=}{\beta_i^{n,k}}{n+k>i}{\beta_i^{n,k}+\beta_{3n+k-i}^{n,k}}
\end{equation}

The same method can be used to compute $\beta_i^{n,k}$ mod $p$ for any prime $p>2$. All we have to show is that $\beta_i^{1,0}=-i$ for $1\le i\le p-1$. The easiest way to see this is to note that $\beta_i^{1,0}=\alpha_{(p,i)}^{1,i}$ for $1\le i\le p-1$, which is the same as the coefficient of $e_pe_1^i$ in the expansion of $m_{(p,i)}=m_pm_i-m_{p+i}$ (by the symmetry of the transition matrix from the basis $\{m_\lambda\}$ to $\{e_\lambda\}$); then a simple application of the Waring formula gives the answer. Everything else follows exactly in the same way as above.
}\end{exa}

\begin{exa}{\rm
Let us compute the coefficient of $e_1^2e_2e_{5n+k-4}$ in $m_{(1^k5^n)}$ for all $n,k$ with $5n+k>4$. 

For $n=1$, we have $\beta_{(1^22)}^{1,0}+4\beta_{(1^4)}^{1,0}+\beta_{(1^32)}^{1,0}=0$, whence $\beta_{(1^22)}^{1,0}=-4$. On the other hand, $\beta_{(1^3)}^{1,0}+\beta_{(1^4)}^{1,0}=0$, whence $\beta_{(1^3)}^{1,0}=-1$. In consequence, 
\[\alpha_{(1^32)}^{1,0}=\beta_{(1^22)}^{1,0}+\beta_{(1^3)}^{1,0}=-5\,,\quad\alpha_{(1^22(k+1))}^{1,k}=\beta_{(1^22)}^{1,0}=-4\,,\spa k\ge 1\,.\]

For $n=2$, we have to compute first $\beta_{(1^22)}^{2,0}+4\beta_{(1^4)}^{2,0}+\beta_{(1^32)}^{2,0}=0$. But $\beta_{(1^4)}^{2,0}=\beta_\emptyset^{1,0}=5$ (the coefficient of $e_5$ in $m_5$), and $\beta_{(1^32)}^{2,0}=\beta_1^{1,0}=-1$, by the observation at the end of Example \ref{exaone}. 
Hence $\beta_{(1^22)}^{2,0}=-19$. Similarly, we obtain $\beta_{(1^26)}^{2,0}=4$. 
In consequence,
\[\alpha_{(1^226)}^{2,0}=\beta_{(1^22)}^{2,0}+\beta_{(1^26)}^{2,0}=-15\,,\quad\alpha_{(1^22(k+6))}^{2,k}=-19-4k\,,\spa k\ge 1\,\,.\]

For $n\ge 3$, we only have to compute $\beta_{(1^22)}^{n,0}=\alpha_{(1^22(5n-4))}^{n,0}$. If we are only interested in working mod 5, then we can apply Corollary \ref{simpleex} and obtain
\[\alpha_{(1^22(5n+k-4))}^{n,k}=\frac{2n+k-3}{n-1}\,\binom{n+k-3}{n-2}\quad \mbox{(mod 5)}\,,\qquad\mbox{for}\spa n\ge 3,\,k\ge 0\,.\]
}\end{exa}

\section{Steenrod Operations in the Schubert Picture of the Cohomology of Grassmannians}

In this section, we investigate the action of the operators $D_{(k^n)}$ in $\du$ (in particular, $SQ^n$) on Schur functions. As we have pointed out in \S3, this corresponds to the action of Steenrod operations on Schubert classes in the cohomology of Grassmannians. More precisely, for every prime $p$, the mod $p$ reduction of the action of $D_{((p-1)^n)}$ gives the action of $\calp^n$; in particular, the mod 2 reduction of the action of $SQ^n$ gives the action of $Sq^{2n}$. 

Let us start by recalling some terminology concerning Young diagrams. The {\em content} $c(b)$ of a box $b$ in row $i(b)$ and column $j(b)$ is the difference $j(b)-i(b)$. The {\em hooklength} $h(b)$ is defined to be $\lambda_{i(b)}+\lambda_{j(b)}'-i(b)-j(b)+1$. A shape $\lambda$ is called a {\em border strip} if it is connected and does not contain a $2\times 2$ block of boxes. A shape satisfying just the second condition is called a {\em broken border strip}; clearly, such a shape is a union of connected components, each of which is a border strip. The {\em height} ht$(\lambda)$ of a broken border strip $\lambda$ is defined to be one less than the number of rows it occupies. A {\em sharp corner} in a broken border strip is a box having no boxes above it or to its left. A {\em dull corner} in a broken border strip is a box that has a box to its left and a box above it, but no box directly northwest of it. We denote by SC$(\lambda)$ and DC$(\lambda)$ the sets of sharp and dull corners of the broken border strip $\lambda$, and by $cc(\lambda)$ the number of its connected components; if the skew shape $\lambda$ is not a broken border strip, it is convenient to set $cc(\lambda):=\infty$. 

\begin{prop}\label{actdk}
We have that
\[D_k\,s_\lambda=\stacksum{\mu\supset\lambda\,:\,|\mu|-|\lambda|=k}{cc(\mu/\lambda)\le 2} d_{\lambda\mu}\,s_\mu\,,\qquad D_k^\perp\,s_\lambda=\stacksum{\mu\subset\lambda\,:\,|\lambda|-|\mu|=k}{cc(\lambda/\mu)\le 2} d_{\mu\lambda}\,s_\mu\,;\]
if $\mu/\lambda$ is a border strip, then 
\[d_{\lambda\mu}=(-1)^{\operatorname{ht}(\mu/\lambda)}\left(\sum_{b\in\operatorname{SC}(\mu/\lambda)}{c}(b)-\sum_{b\in\operatorname{DC}(\mu/\lambda)}{c}(b)\right)\,,\]
and if $\mu/\lambda$ has two connected components, then $d_{\lambda\mu}=(-1)^{\operatorname{ht}(\mu/\lambda)}$.
\end{prop}
\begin{proof} 
We start by recalling the following formula (see \cite{macsfh}):
\begin{equation}\label{macone}
p_k \,s_\lambda=\stacksum{\mu\supset\lambda\,:\,|\mu|-|\lambda|=k}{cc(\mu/\lambda)=1}(-1)^{\operatorname{ht}(\mu/\lambda)}\,s_\mu\,.
\end{equation}
By duality, we have a similar formula:
\begin{equation}\label{mactwo}
p_k^\perp \,s_\lambda=\stacksum{\mu\subset\lambda\,:\,|\lambda|-|\mu|=k}{cc(\lambda/\mu)=1}(-1)^{\operatorname{ht}(\lambda/\mu)}\,s_\mu\,.
\end{equation}
We now combine these formulas with the expression for $D_k$ in Proposition \ref{expd}. We deduce that $D_k\,s_\lambda$ is a linear combination of $s_\mu$, where three cases are possible: (1) both $\pi:=\mu/\lambda$ and $\rho:=\lambda/\mu$ are border strips; (2) $\lambda$ is contained in $\mu$ and $\pi:=\mu/\lambda$ is a border strip; (3) $\lambda$ is contained in $\mu$ and $\pi:=\mu/\lambda$ is a broken border strip with two connected components. These cases are illustrated in the three figures below, where the black boxes are removed, the bold ones are added, and the shaded ones are first removed (by applying the operator $p_i^\perp$), and then added (by applying the operator $p_{k+i}$).

\[
\begin{array}{lcr}
\mbox{\psfig{file=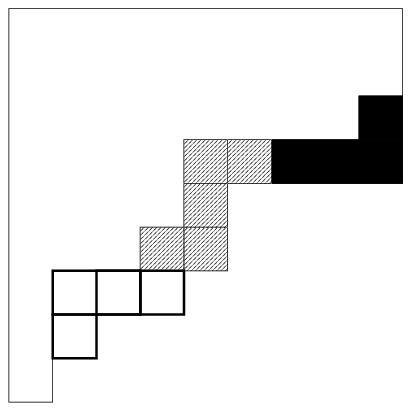}}\spa\spa&\spa\spa\mbox{\psfig{file=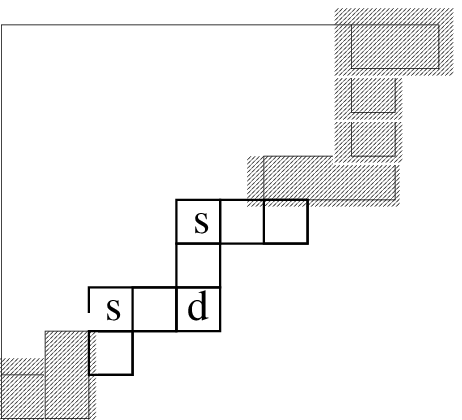}}\spa\spa&\spa\spa\mbox{\psfig{file=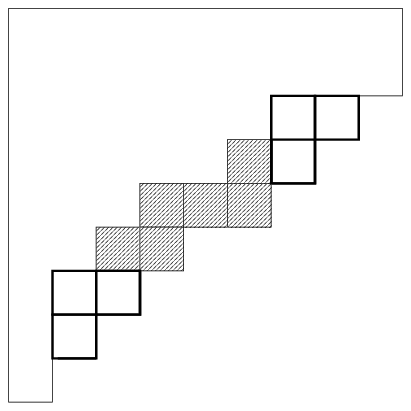}}
\end{array}
\]

%
%
In the first case, let $\theta\subset \lambda$ be the border strip ``connecting'' $\pi$ and $\rho$, i.e. the border strip characterized by the fact that $\pi\cup\theta\cup\rho$ is connected (the shaded area in the first figure). The partition $\mu$ can be obtained in two ways from $\lambda$: by removing $\rho$ and adding $\pi$, or by removing $\rho\cup\theta$ and adding $\pi\cup\theta$. Note that $\mbox{ht}(\pi\cup\theta)=\mbox{ht}(\pi)+\mbox{ht}(\theta)+1$ and $\mbox{ht}(\rho\cup\theta)=\mbox{ht}(\rho)+\mbox{ht}(\theta)$ if $\pi$ is southwest of $\rho$, and viceversa, otherwise. Hence, the corresponding coefficients given by the two formulas above cancel.

In the second case, the partition $\mu$ can be obtained in several ways, by removing a border strip $\theta$ for which $\pi\cup\theta$ is connected, and then adding $\pi\cup\theta$. The border strip $\theta$ can be situated northeast of $\pi$, in which case $\mbox{ht}(\pi\cup\theta)=\mbox{ht}(\pi)+\mbox{ht}(\theta)+1$, or southwest of $\pi$, in which case $\mbox{ht}(\pi\cup\theta)=\mbox{ht}(\pi)+\mbox{ht}(\theta)$. In the first case, the sum of the corresponding coefficients is $(-1)^{\operatorname{ht}(\pi)+1}(i(b_1)-1)$, while in the second case, the sum is $(-1)^{\operatorname{ht}(\pi)}(j(b_2)-1)$; here $b_1$ and $b_2$ are the first and last sharp corners of $\pi$ (from northeast to southwest). On the other hand, it is easy to see that  
\begin{align*}
\sum_{b\in\operatorname{SC}(\pi)}{c}(b)-\sum_{b\in\operatorname{DC}(\pi)}{c}(b)&=\left(\sum_{b\in\operatorname{SC}(\pi)}j(b)-\sum_{b\in\operatorname{DC}(\pi)}j(b)\right)\\
&-\left(\sum_{b\in\operatorname{SC}(\pi)}i(b)-\sum_{b\in\operatorname{DC}(\pi)}i(b)\right)=j(b_2)-i(b_1)\,.\end{align*}

In the third case, there is a unique way to obtain $\mu$ from $\lambda$, namely by removing the border strip $\theta$ ``connecting'' $\pi$, and then adding $\pi\cup\theta$. Since $\mbox{ht}(\pi\cup\theta)=\mbox{ht}(\pi)+\mbox{ht}(\theta)$, we have $d_{\lambda\mu}=(-1)^{\operatorname{ht}(\pi)}$, as claimed.
\end{proof}

\begin{rema}{\rm
Proposition \ref{actdk} has applications to the inverse Kostka matrix $K^{(-1)}$ (see \cite{macsfh}). For example, we can immediately deduce the following closed formula (appearing in \cite{macsfh} p. 107) for $K_{\lambda\mu}^{(-1)}$ with $\lambda=(1^kn)$, $n>1$:
\[\caset{K_{\lambda\mu}^{(-1)}}{=}{(-1)^{n+1}\,(k+1)}{\mu=(1^{k+n})}{(-1)^{n-l}}{\mu=(1^i 2^jl)\,,\spa j+l\le n}{0}\]
indeed, we only need to notice that $D_{n-1}\,e_{k+1}=m_{(1^kn)}$, which means that in this case $K_{\lambda\mu}^{(-1)}=d_{(1^{k+1})\mu}$ for all $\mu$. In fact, we can use this idea more generally, and give an efficient algorithm for computing the inverse Kostka matrix. For every partition $\lambda$, we define an operator $\cald_\lambda$ in $\du$ recursively as follows:
\[\cald_k:=D_{k-1}\,,\quad \cald_\lambda:=D_{\lambda_1-1}\,\cald_{\lambda\setminus\lambda_1}-\sum_{i=2}^{l(\lambda)}\lambda_i\,\cald_{(\lambda\setminus\lambda_1\setminus\lambda_i)\cup(\lambda_1+\lambda_i-1)}\spa\mbox{for } l(\lambda)\ge 2\,.\]
It is easy to see by induction that $\cald_\lambda\,e_{l(\lambda)}=m_\lambda$. On the other hand, we can express $\cald_\lambda\,e_{l(\lambda)}$ in the basis of Schur functions using the formula in Proposition \ref{actdk} repeatedly. Thus, we determine $K_{\lambda\mu}^{(-1)}$ for all partitions $\mu$.
}\end{rema}

In order to express the action of $D_{(k^{n})}$, we need the following identity in $\du$:
\begin{equation}\label{basic}
n\,D_{(k^{n})}=\sum_{i=1}^n (-k-1)^{i-1}\,D_{ik}\circ D_{(k^{n-i})}\,.
\end{equation}
R. Wood proves this identity in \cite{woodos} for $k=1$, using his composition formula for two basis elements in $\du$; his argument can be extended without difficulty to any $k$. By combining this recursive formula for $D_{(k^n)}$ with Proposition \ref{actdk}, we obtain the following Corollary.

\begin{cor}\label{actdkn}
We have that
\[D_{(k^n)}\,s_\lambda=\sum_{\mu\supset\lambda\,:\,|\mu|-|\lambda|=kn} a_{\lambda\mu}\,s_\mu\,;\]
the coefficients $a_{\lambda\mu}$ are given by
\[a_{\lambda\mu}=\sum_{\lambda=\mu_{(0)}\subset\mu_{(1)}\subset\ldots\subset\mu_{(m)}=\mu}f(\mu_{(0)},\mu_{(1)})\,f(\mu_{(1)},\mu_{(2)})\,\ldots\,f(\mu_{(m-1)},\mu_{(m)})\,;\]
here $\mu_{(i)}/\mu_{(i-1)}$ is a broken border strip with at most two connected components whose length is a multiple of $k$ for every $i$, and
\[f(\pi,\rho)=\frac{k\,(-k-1)^{(|\rho|-|\pi|)/k-1}}{|\rho|-|\lambda|}\,d_{\pi\rho}\,.\]
\end{cor}

\begin{rema}{\rm
Since $Sq^{2n}$ squares every cohomology class of dimension $2n$, the coefficient $a_{\lambda\mu}$ for $k=1$ and $|\mu|=2|\lambda|$ is congruent mod 2 to the Littlewood-Richardson coefficient $c_{\lambda\lambda}^\mu$. Thus all our results for $a_{\lambda\mu}$ (Proposition \ref{actdkn}, Theorem \ref{main}, Corollaries \ref{corone} and \ref{cortwo}) can be translated in terms of mod 2 congruences for $c_{\lambda\lambda}^\mu$.
}\end{rema}

Corrolary \ref{actdkn} has an important consequence concerning the way in which the coefficients $a_{\lambda\mu}$ depend on the embedding of the skew diagram $\mu/\lambda$ and its connected components in the plane. In order to state this fact, we need additional notation. Let $\pi=\pi(\lambda,\mu)$ and $\rho=\rho(\lambda,\mu)$ be the diagrams of smallest weight ($\pi$ possibly empty) for which $\mu/\lambda=\rho/\pi$, and let $c(\lambda,\mu)$ be the content of the top left box of $\rho(\lambda,\mu)$ (in its embedding in $\mu$). Similarly, we write the connected components of $\mu/\lambda$ (from northeast to southwest say) as $\rho_{(i)}(\lambda,\mu)/\pi_{(i)}(\lambda,\mu)$, $i=1,2,\ldots$, where the corresponding Young diagrams are again of minimum weight. The corresponding contents are denoted by $c_i(\lambda,\mu)$. 

\begin{cor}
The coefficient $a_{\lambda\mu}$ depends only on the connected components of $\mu/\lambda$ and the corresponding contents $c_i(\lambda,\mu)$. In particular, it depends only on $\mu/\lambda$ and $c(\lambda,\mu)$. 
\end{cor}
\begin{proof}
According to Proposition \ref{actdk}, the coefficient $d_{\lambda\mu}$ depends on $c(\lambda,\mu)$ if and only if $\mu/\lambda$ is connected. The Corollary follows by induction on $|\mu|-|\lambda|$ by combining the previous observation with (\ref{basic}).
\end{proof}

From now on we concentrate on the case $k=1$. In other words we work with the operators $SQ^n=D_{(1^n)}$, which reduce to the Steenrod squares $Sq^{2n}$ upon reduction mod 2. In order to investigate the way in which $a_{\lambda\mu}$ depends on $c(\lambda,\mu)$, we need the following two lemmas.

\begin{lem}\label{lllooo}
The operator $\sum_{i\ge 1} p_i\,p_i^\perp$ multiplies every homogeneous symmetric function by its degree. In particular, we have the following generalization of Newton's identity (which can be recovered by choosing $\lambda=(1^n)$):
\[|\lambda|\,s_\lambda=\sum_{\mu\subset\lambda\,:\,cc(\lambda/\mu)=1}(-1)^{\operatorname{ht}(\lambda/\mu)}\,p_{|\lambda|-|\mu|}\,s_\mu\,.\]
\end{lem}
\begin{proof}
Using the same technique as in the proof of Proposition \ref{expd}, we can show that the operator $\sum_{i\ge 1} p_i\,p_i^\perp$ coincides with the operator $\sum_{i\ge 1} x_i\,\partial_i$. We can easily check that the latter multiplies every monomial in the $x_i$'s by its degree. The identity follows by combining this result with identity (\ref{mactwo}).
\end{proof}

\begin{lem}\label{lllttt}
For every symmetric functions $f,g$, and every $k>0$, we have
\[f^\perp\,(p_kg)=p_k\,(f^\perp\,g)+(p_k^\perp\,f)^\perp\,g\,.\]
In particular, we have the following generalization of (\ref{macone}):
\[p_k\,s_{\mu/\lambda}=\stacksum{\nu\supset\mu\,:\,|\nu|-|\mu|=k}{cc(\nu/\mu)=1}(-1)^{\operatorname{ht}(\nu/\mu)}\,s_{\nu/\lambda}-\stacksum{\nu\subset\lambda\,:\,|\lambda|-|\nu|=k}{cc(\lambda/\nu)=1}(-1)^{\operatorname{ht}(\lambda/\nu)}\,s_{\mu/\nu}\,.\]
\end{lem}
\begin{proof}
The first statement holds more generally. Let $B$ be a bialgebra with a countable basis $\{c_\omega\}$ containing 1, and $B^\circ$ its finite dual with basis $\{c^\omega\}$, such that $\dual{c^\pi}{c_\rho}=\delta_{\pi\rho}$. Consider elements $f\in B^\circ$ and $g,p\in B$, such that $p$ is primitive. Recall the canonical left action of $B^\circ$ on $B$ given by $f\rightharpoonup g:=\sum \dual{f}{g_2}\,g_1$, where we use the Sweedler notation for the comultiplication. Similarly, there is a canonical left action of $B$ on $B^\circ$, for which the same notation is traditional. The action of symmetric functions on themselves via Hammond operators are special cases of the canonical actions of bialgebras on their duals. Since $p$ is primitive, we have
\[\case{c^\omega\rightharpoonup p}{=}{p}{c^\omega=1}{\dual{c^\omega}{p}}\]
Now recall the fact that $B$ is a $B^\circ$-module algebra under $\rightharpoonup$, and that the comultiplication can be expressed in terms of the basis $\{c_\omega\}$ as follows:
\[\delta(f)=\sum_\omega (c_\omega\rightharpoonup f)\otimes c^\omega\,.\]
Combining these facts, we have
\begin{align*}
f\rightharpoonup (gp)-(f\rightharpoonup g)p&=\sum_{\omega\,:\,c_\omega\ne 1} \dual{c^\omega}{p}\,((c_\omega\rightharpoonup f)\rightharpoonup g)\\
&=\left(\left(\sum_\omega\dual{c^\omega}{p}\,c_\omega\right)\rightharpoonup f\right)\rightharpoonup g=(p\rightharpoonup f)\rightharpoonup g\,.
\end{align*}

The formula for $p_k\,s_{\mu/\lambda}$ then follows easily by letting $f=s_\lambda$, $g=s_\mu$, and by combining the previous result with (\ref{macone}) and (\ref{mactwo}).
\end{proof}

Let us now recall the principal specialization (in fact, a special case of it) of Schur functions. For every non-negative integer $n$, we can define a map $\operatorname{ps}_n^1$ from $Sym_*$ to $\bZ$ which sends a symmetric function $f(x_1,x_2,\ldots)$ to the integer obtained by setting $x_1=\ldots=x_n=1$ and $x_k=0$ for $k>n$. In fact, we can define $\operatorname{ps}_n^1$ for every integer $n$, by letting $\operatorname{ps}_{n}^1(f)=\operatorname{ps}_{-n}^1(S(f))$ if $n<0$, where $S$ is the antipode in $Sym_*$. It is easy to see that $\operatorname{ps}_n^1$ are ring homomorphisms, and that $\operatorname{ps}_n^1(p_k)=n$ for every $k$. It is also known (see \cite{macsfh}) that
\begin{equation}\label{spec}
\operatorname{ps}_n^1(s_{\lambda'})=\prod_{b\in\lambda}\frac{n-{c}(b)}{h(b)}\,.
\end{equation}

We are now able to prove our result concerning the way in which $a_{\lambda\mu}$ depends on $c(\lambda,\mu)$. 

\begin{thm}\label{main}
We have
\[a_{\lambda\mu}=\operatorname{ps}_{c(\lambda,\mu)}^1\left((\cals^{|\pi(\lambda,\mu)|}\,s_{\pi'(\lambda,\mu)})^\perp\,s_{\rho'(\lambda,\mu)}\right)\,,\]
where $\cals^n=\sum_{i=0}^n (-1)^i\,SQ^i$.
\end{thm}
\begin{proof}
We let $\pi=\pi(\lambda,\mu)$, $\rho=\rho(\lambda,\mu)$, $c=c(\lambda,\mu)$, and $n=|\rho|-|\pi|$, for simplicity. We apply induction on $n$, which starts at $n=1$ by Proposition \ref{actdk}. Identity (\ref{basic}) tells us that we have
\begin{equation}\label{idone}
na_{\lambda\mu}=\sum_{\lambda\subseteq\nub\subset\mu\,:\,cc(\mu/\nub)\le 2}\,(-2)^{|\mu|-|\nub|-1}a_{\lambda\nub}\,d_{\nub\mu}\,.
\end{equation}
If $\mu/\nub$ is a border strip, then $d_{\nub\mu}=(-1)^{\operatorname{ht}(\rho/\nu)}c+d_{\nu\rho}$, where $\nu:=\nub/\pi$; otherwise $d_{\nub\mu}=d_{\nu\rho}$. We want to lift (\ref{idone}) to an identity for symmetric functions. The crucial ingredients for our proof are the symmetric functions 
\[q_k:=\frac{p_k+p_{k-1}+2p_{k-2}+\ldots+2^{k-2}p_1}{2^{k-1}}\,,\]
which clearly satisfy $\operatorname{ps}_c^1(q_k)=c$. We claim that the following identity holds in $Sym_*\otimes\bQ$:
\begin{align}\label{idtwo}
n\,(\cals^{|\pi|}\,s_{\pi'})^\perp\,s_{\rho'}&=\stacksum{\pi\subseteq\nu\subset\rho}{cc(\rho/\nu)=1}\,(-2)^{|\rho|-|\nu|-1}\,\left((-1)^{\hgt(\rho/\nu)}q_{|\rho|-|\nu|}+d_{\nu\rho}\right)\,(\cals^{|\pi|}\,s_{\pi'})^\perp\,s_{\nu'}\\
&+\stacksum{\pi\subseteq\nu\subset\rho}{cc(\rho/\nu)=2}\,(-2)^{|\rho|-|\nu|-1}\,d_{\nu\rho}\,(\cals^{|\pi|}\,s_{\pi'})^\perp\,s_{\nu'}\,.\nonumber
\end{align}
If we prove this identity, the induction step is straightforward: we simply apply $\operatorname{ps}_c^1$ to it, and compare the result with (\ref{idone}).

Identity (\ref{idtwo}) is an identity for non-homogeneous symmetric functions, so it breaks down into identities in every degree. In degree $n-i$ ($0\le i\le n$), we have
\begin{gather}\label{idthree}
(-1)^in\,(SQ^i\,s_{\pi'})^\perp\,s_{\rho'}=\sum_{j=0}^i\stacksum{\nu\subset\rho}{cc(\rho/\nu)=1}\,(-1)^{j-1+|\rho|-|\nu|+\hgt(\rho/\nu)}\,2^{i-\alpha(j)-1}\,p_{|\rho|-|\nu|+j-i}\\
\times(SQ^j\,s_{\pi'})^\perp\,s_{\nu'}+\sum_{j=0}^{i-1}\stacksum{\nu\subset\rho\,:\,|\nu|=|\rho|+j-i}{cc(\rho/\nu)\le 2}\,(-1)^{i-1}\,2^{i-j-1}\,d_{\nu\rho}\,(SQ^j\,s_{\pi'})^\perp\,s_{\nu'}\,;\nonumber
\end{gather}
here $\alpha(j)$ is $j$ or $j-1$ depending on $j$ being less than $i$ or equal to $i$, respectively, and we adopt the convention that $p_n=0$ if $n\le 0$. Note that we dropped the restriction $\pi\subseteq\nu$, because if $\nu$ does not satisfy it, then $(SQ^j\,s_{\pi'})^\perp\,s_{\nu'}=0$ (recall that the partitions indexing the Schur functions in the expansion of $SQ^j\,s_{\pi'}$ contain $\pi'$). For every $j=0,\ldots,i-1$, we pair the corresponding terms in two sums in the RHS of (\ref{idthree}). We start by investigating the $j$-th term in the second sum. We have
\[\stacksum{\nu\subset\rho\,:\,|\nu|=|\rho|+j-i}{cc(\rho/\nu)\le 2}d_{\nu\rho}\,s_\nu=D_{i-j}^\perp\,s_\rho\,.\]
We express the action of $D_{i-j}^\perp$ on $s_\rho$ using Proposition \ref{expd} and (\ref{mactwo}), then we apply the standard involution on symmetric functions to the result, and obtain
\[\stacksum{\nu\subset\rho\,:\,|\nu|=|\rho|+j-i}{cc(\rho/\nu)\le 2}d_{\nu\rho}\,s_{\nu'}=\stacksum{\nu\subset\rho}{cc(\rho/\nu)=1}(-1)^{j-i-1+|\rho|-|\nu|+\hgt(\rho/\nu)}\,p_{|\rho|-|\nu|+j-i}\,s_{\nu'}\,.\]
By Lemma \ref{lllttt}, when we apply the Hammond operator $(SQ^j s_{\pi'})^\perp$ to the previous expression, we can express the result as a sum of two terms, one of which is the $j$-th term in the first sum of the RHS of (\ref{idthree}) up to a constant. Hence, by pairing the $j$-th terms in the two sums, we obtain
\[(-1)^{j+|\rho|}\,2^{i-j-1}\,\stacksum{\nu\subset\rho}{cc(\rho/\nu)=1}\,(-1)^{\hgt(\rho/\nu)-|\nu|}\,\left(p_{|\rho|-|\nu|+j-i}^\perp\,(SQ^j\,s_{\pi'})\right)^\perp\,s_{\nu'}\,.\]
We can simplify this sum by using (\ref{mactwo}) again, together with the standard involution on symmetric functions:
\begin{equation}\label{idfour}
\stacksum{\nu\subset\rho\,:\,|\nu|=k}{cc(\rho/\nu)=1}(-1)^{\hgt(\rho/\nu)}s_{\nu'}=(-1)^{|\rho|-k-1}\,p_{|\rho|-k}^\perp\,s_{\rho'}\,.
\end{equation}
Hence, by the commutativity of Hammond operators, the sum of the $j$-th terms in the two sums in the RHS of (\ref{idthree}) can be expressed as
\begin{equation}\label{idzz}
(-1)^{j-1}\,2^{i-j-1}\left(\left(\sum_k p_{|\rho|-k}\,p_{|\rho|-k+j-i}^\perp\right)\,(SQ^j\,s_{\pi'})\right)^\perp s_{\rho'}\,;
\end{equation}
here the summation ranges over all $k$ with $|\rho|-k$ lying between 1 and the length of the longest border strip of $\rho$. We want to show that we can actually extend the range of $|\rho|-k$ to infinity without affecting the result; this has the advantage of being able to write the sum as the operator $D_{i-j}$. We regard $SQ^j\,s_{\pi'}$ as a linear combination of Schur functions $s_\sigma$ with  $|\sigma|=|\pi|+j$ and $\pi'\subseteq\sigma$. Applying $D_{i-j}$ viewed as $\sum_{l\ge 1}p_l\,p_{l+j-i}^\perp$ to $s_\sigma$ gives a sum of Schur functions $s_\tau$, where $\tau$ is obtained from $\sigma$ by removing a border strip of length $l+j-i>0$, and then adding a border strip of length $l$. Since we then consider $s_{\rho'/\tau}$, we are only interested in those $\tau$ contained in $\rho'$. Clearly, for such $\tau$, the corresponding added border strip has length $l$ at most the length of the longest border strip of $\rho$, that is within the range of $|\rho|-k$. According to the above remarks and (\ref{basic}), we can express the RHS of (\ref{idthree}) less the term corresponding to $j=i$ in the first sum as
\begin{equation}\label{idfive}
\sum_{j=0}^{i-1} (-1)^{j-1}\,2^{i-j-1}\,((D_{i-j}\,SQ^j)\,s_{\pi'})^\perp\,s_{\rho'}=(-1)^i\,i\,(SQ^i\,s_{\pi'})^\perp\,s_{\rho'}\,.
\end{equation}
It only remains to deal with the term corresponding to $j=i$ in the RHS of (\ref{idthree}). Applying (\ref{idthree}) again, the commutativity of Hammond operators, and Lemma \ref{lllooo}, we can express this term as follows:
\begin{gather}\label{idsix}
(-1)^i \,\stacksum{\nu\subset\rho}{cc(\rho/\nu)=1}\,(-1)^{|\rho|-|\nu|-1+\hgt(\rho/\nu)}\,p_{|\rho|-|\nu|}\,(SQ^i\,s_{\pi'})^\perp\,s_{\nu'}\\
=(-1)^i\left(\sum_k p_{|\rho|-k}\,p_{|\rho|-k}^\perp\right)\left((SQ^i\,s_{\pi'})^\perp\,s_{\rho'}\right)=(-1)^i(n-i)\, (SQ^i\,s_{\pi'})^\perp\,s_{\rho'}\,;\nonumber
\end{gather}
here the range of $|\rho|-k$ is the same as in (\ref{idzz}), and was extended to infinity by the same argument as above. Combining (\ref{idfive}) and (\ref{idsix}) proves (\ref{idthree}).
\end{proof}

Theorem \ref{main} says that $a_{\lambda\mu}$ is a certain specialization of a non-homogeneous symmetric function with highest homogeneous component $s_{\mu'/\lambda'}$. More precisely, $a_{\lambda\mu}$ is a polynomial in $c(\lambda,\mu)$ depending only on $\mu/\lambda$; it has degree $|\mu|-|\lambda|$, free term equal to $a_{\pi(\lambda,\mu)\rho(\lambda\mu)}$, and coefficient of the leading term equal to the corresponding coefficient of $\operatorname{ps}_{c(\lambda,\mu)}^1(s_{\mu'/\lambda'})$. Thus, we have reduced the computation of $a_{\lambda\mu}$ to the case $l(\lambda)<l(\mu)$, $\lambda_1<\mu_1$. One can use Proposition \ref{actdkn} to compute $SQ^i\,s_{\pi'(\lambda,\mu)}$ for all $i=1,\ldots,|\pi(\lambda,\mu)|$ (in fact we only need the coefficients $a_{\pi'(\lambda,\mu)\nu}$ with $\nu\subseteq\rho'(\lambda,\mu)$). Then Theorem \ref{main} enables us to compute quite easily the coefficient $a_{\lambda\mu}$. Of course, the complexity is in the previous step, and it increases with the weight of $\pi(\lambda,\mu)$. The following Corollary considers the cases when this weight is 0 or 1.

\begin{cor}\label{corone}
If $\nu:=\mu/\lambda$ is a Young diagram, then $a_{\lambda\mu}=\operatorname{ps}_{c(\lambda,\mu)}^1(s_{\nu'})$, which can be expressed using (\ref{spec}). If $\pi(\lambda,\mu)$ consists of a single box, then $a_{\lambda\mu}=\operatorname{ps}_{c(\lambda,\mu)}^1((p_1-p_2)^\perp\,s_{\rho'(\lambda,\mu)})$.
\end{cor}
\begin{proof}
The above statements are immediate consequences of Theorem \ref{main}. All we need to observe (for the second part) is that $SQ^1\,s_1=D^1\,p_1=p_2$, by Proposition \ref{expd}.
\end{proof}

The following Corollary gives sufficient conditions for a coefficient $a_{\lambda\mu}$ to be 0. To state it, we need additional notation. For every partition $\lambda$ and non-negative integer $n$, we denote by $m(\lambda,n)$ the largest number $k$ for which either $\lambda_{k}'+\ldots+\lambda_{l(\lambda)}'>n$ or $\lambda_{k+1}'+\ldots+\lambda_{l(\lambda)}'=n$; if no such $k$ exists, we let $m(\lambda,n):=0$. 

\begin{cor}\label{cortwo}
Given diagrams $\lambda\subset\mu$, we consider the partitions $\sigma$ and $\tau$ whose parts are the lengths of columns and rows of $\mu/\lambda$, respectively. The coefficient $a_{\lambda\mu}$ is 0 whenever the following condition is satisfied:
\[1-m(\sigma,|\pi(\lambda,\mu)|)\le c(\lambda,\mu)\le m(\tau,|\pi(\lambda,\mu)|)-1\,.\]
\end{cor}
\begin{proof}
We let $\pi=\pi(\lambda,\mu)$ and $\rho=\rho(\lambda,\mu)$, for simplicity. Clearly, it suffices to consider the case $2|\pi|<|\rho|$, since otherwise $m(\sigma,|\pi|)=m(\theta,|\pi|)=0$, so the above condition cannot hold. Given $i$ between 0 and $|\pi|$, we view $(SQ^i\,s_{\pi'})^\perp\,s_{\rho'}$ as a linear combination of skew Schur functions $s_{\rho'/\nu}$, where $\pi'\subseteq\nu\subseteq\rho'$ and $|\nu|-|\pi|=i$. We then view every $s_{\rho'/\nu}$ as a linear combination (with non-zero coefficients) of Schur functions $s_{\theta'}$. By the Littlewood-Richardson rule, the minimum value $m_c(\nu)$ of $\theta_1$ for a fixed $\nu$ is equal to the length of the longest column of the skew shape $\rho'/\nu$. Indeed, this is clearly a lower bound, since the boxes in the mentioned column have to be filled with different symbols in any Littlewood-Richardson filling; on the other hand, the lower bound is attained, since by filling every column of $\rho'/\nu$ with symbols $1,2,\ldots$ (from top to bottom, in this order), we obtain a Littlewood-Richardson filling. By applying the standard involution on symmetric functions, we can deduce from here that the minimum value $m_r(\nu)$ of $l(\theta)$ for a fixed $\nu$ is equal to the length of the longest row of the skew shape $\rho'/\nu$. By (\ref{spec}), we have that $\operatorname{ps}_{c(\lambda,\mu)}^1(s_{\rho'/\nu})=0$ whenever $1-m_r(\nu)\le c(\lambda,\mu)\le m_c(\nu)-1$. To minimize $m_r(\nu)$, we let $i=|\pi|$, and construct a sequence of Young diagrams $\pi'=\nu_{(0)},\nu_{(1)},\ldots,\nu_{(|\pi|)}$ such that every $\nu_{(j)}$ is obtained by adding a box to one of the rows of $\rho'/\nu_{(j-1)}$ of maximum length. It is not difficult to see that there is always such a choice, and that $m_r(\nu_{(|\pi|)})=m(\sigma,|\pi|)$ minimizes $m_r(\nu)$. The same argument applies to $m_c(\nu)$, thus proving the Corollary.
\end{proof}

\[
\begin{array}{c}
\mbox{\psfig{file=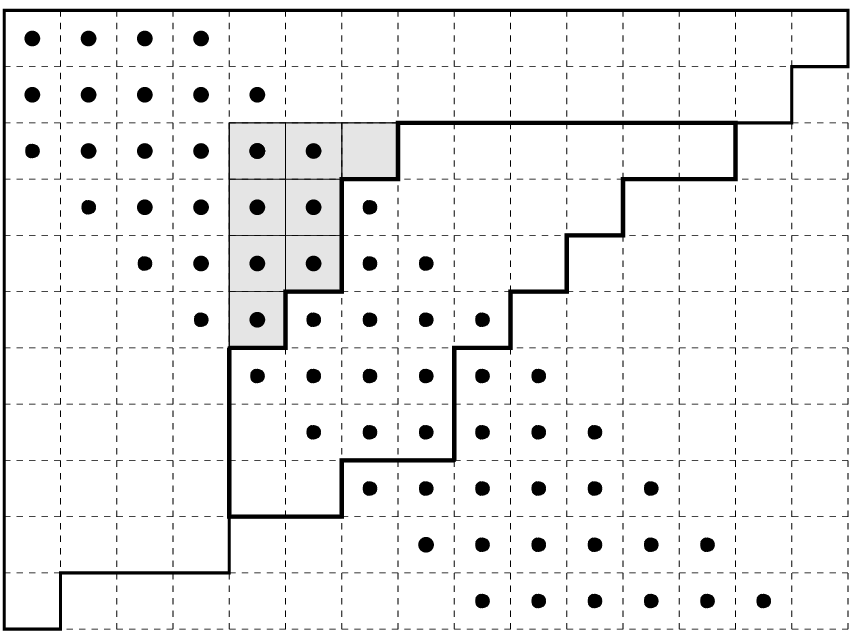}}
\end{array}
\]

Note that if $|\mu|-|\lambda|\le 2$, the condition in the Corollary is also necessary (which is not true in general, however). Indeed, if $|\mu|-|\lambda|=1$, the necessity is clear. If $|\mu|-|\lambda|=2$, there are three cases possible: (1) $\mu/\lambda=(2)$, (2) $\mu/\lambda=(1,1)$, and (3) $\mu/\lambda$ is a skew diagram with two boxes of content $c_1,c_2$. By Proposition \ref{actdkn}, we have $a_{\lambda\mu}=c(\lambda,\mu)\,(c(\lambda,\mu)-1)/2$ in the first case, $a_{\lambda\mu}=c(\lambda,\mu)\,(c(\lambda,\mu)+1)/2$ in the second one, and $a_{\lambda\mu}=c_1c_2+1$ in the third one. The necessity of the condition in the Corollary follows without difficulty from here. 

\begin{exa}{\rm
Let $\pi=(3,2,2,1)$, $\rho=(9,7,6,5,4,4,2)$, and $\lambda,\mu$ be such that $\pi(\lambda,\mu)=\pi$ and $\rho(\lambda,\mu)=\rho$. We have $\sigma=(6,5,4,4,3,3,2,1,1)$ and $\tau=(6,5,4,4,4,4,2)$, whence $m(\sigma,8)=3$ and $m(\tau,8)=4$. According to the Corollary, we have $a_{\lambda\mu}=0$ for all $\lambda,\mu$ with $-2\le c(\lambda,\mu)\le 3$. In the figure above, illustrating this example, we marked with a dot the boxes where the top left box of $\rho$ can be translated, such that the corresponding coefficient $a_{\lambda\mu}$ is 0.
}\end{exa}

The next logical step after Theorem \ref{main} would be to investigate how $a_{\lambda\mu}$ depends on the contents $c_i(\lambda,\mu)$ corresponding to the connected components of $\mu/\lambda$ (assuming that this is not connected). It is known that if $\mu/\lambda$ consists of at least two connected components (each of which is a skew diagram), then $s_{\mu/\lambda}$ is the product of the corresponding skew Schur functions. This property does not seem to have a straightforward analogue for the coefficients $a_{\lambda\mu}$, which not only depend on $\mu/\lambda$, but also on $c(\lambda,\mu)$. We conclude this section with a conjecture which represents such an analogue when $\mu/\lambda$ has two connected components. Before stating it, we introduce the following notation involving two Young diagrams $\lambda,\mu$, and two integers $k,n$:
\[f(\lambda,\mu,k,n):=\sum_{\nu\subseteq\lambda,\,\nu\subseteq\mu}\left(\prod_{b\in\lambda/\nu}\frac{k-c(b)}{h(b)}\right)\left(\prod_{b\in\mu/\nu}\frac{n-c(b)}{h(b)}\right)\,;\]
here $\nu$ can be the empty diagram, and if $\nu=\lambda$ (or $\nu=\mu$), the corresponding product is considered to be 1. 

\begin{conj}\label{conj}
Let $\lambda\subset\mu$ be such that $\mu/\lambda$ is a union of two connected components. Let
\[\left(\cals^{|\pi_{(k)}(\lambda,\mu)|}\,s_{\pi_{(k)}'(\lambda,\mu)}\right)^\perp\,s_{\rho_{(k)}'(\lambda,\mu)}=\sum_{\nu\in P_k}\alpha_{k,\nu}\,s_{\nu'}\,,\]
for $k=1,2$, where $P_k$ are corresponding sets of partitions. Then we have
\[a_{\lambda\mu}=\stacksum{\nu_{(1)}\in P_1}{\nu_{(2)}\in P_2} \alpha_{1,\nu_{(1)}}\,\alpha_{2,\nu_{(2)}}\,f(\nu_{(1)},\nu_{(2)},c_1(\lambda,\mu),c_2(\lambda,\mu))\,.\]
\end{conj}

This conjecture was tested for small diagrams $\pi_{(k)},\rho_{(k)}$, $k=1,2$. There are indications that similar formulas exist for three or more connected components. Such formulas would reduce the computation of $a_{\lambda\mu}$ to the case $l(\lambda)<l(\mu)$, $\lambda_1<\mu_1$, and $\mu/\lambda$ having only one connected component, whence they would increase the efficiency of the computation with respect to the formula in Theorem \ref{main}. 

\section{Applications to the Geometry of Grassmannians}

Corollary \ref{corone} provides a closed formula for the action of Steenrod operations on Schubert classes of the form $\sigma_{(j^{n-i}k^i)}$ and $\sigma_{(j^{n-i-1}(j+1)k^i)}$ in the cohomology of $Gr_n(\bC^{n+k})$. Indeed, since Schubert classes in the cohomology of this Grassmannian are indexed by partitions with at most $n$ rows and $k$ columns, all the coefficients $a_{\lambda\mu}$ we need to compute are of the form specified by the Corollary. In general, we can combine the above results (Proposition \ref{actdkn}, Theorem \ref{main}, and possibly Conjecture \ref{conj}), to compute the action of Steenrod operations on any Schubert class. 

Let us consider some examples. Since the cell structure and the action of the Steenrod algebra on projective spaces is well-known, we start with the smallest example which is not a projective space, namely $Gr_2(\bC^4)$ of dimension 8. Here we consider the Schubert cells
\[\varOmega_{1}:=\setp{V}{\bC^1\subset V\subset \bC^3}\,,\quad\varOmega_{(2,1)}:=\setp{V}{\mbox{dim}(V\cap \bC^2) \ge 1}\]
of dimensions 2 and 6, respectively. According to our results, we have $Sq^2\, \sigma_{(2,1)}=0$, which suggests that the top cell $\varOmega_{(2,2)}$ of $Gr_2(\bC^2)$ (the space
itself) is attached trivially to $\varOmega_{(2,1)}$ (see the discussion in \S3). On the other hand, we have
\[Sq^6\, \sigma_1=0\,,\quad (Sq^2 \,Sq^4)\, \sigma_1=0\,,\quad (Sq^2)^3 \,\sigma_1=0\,.\]
The integral lifts of these also give 0, but 
$(SQ^4\circ SQ^2) \,s_1=-2 s_{(2,2)}$ (indeed, the Adem relations do not hold
integrally). There are no higher Steenrod operations to take us from
dimension 2 to 8, but $D_3\,s_1=0$. This suggests that the top cell is attached
trivially to $\varOmega_1$.

Let us now consider $Gr_2(\bC^5)$ of dimension 12, and the Schubert cells
\[\varOmega_{(1,1)}:=\setp{V\subset\bC^3}{\mbox{dim}(V\cap \bC^2) \ge 1}\,,\quad
\varOmega_{(3,1)}:=\setp{V}{\mbox{dim}( V\cap \bC^2)\ge 1}\]
of dimensions 4 and 8, respectively. We have
\[Sq^4\, \sigma_{(3,1)}=0\,,\quad (Sq^2)^2\, \sigma_{(3,1)}=0\,,\quad \calp^1 \,\sigma_{(3,1)}=0\spa (\mbox{for } p=3).\]
The integral lifts of these operations also give 0. This
suggests that the top cell $\varOmega_{(3,3)}$ of $Gr_2(\bC^5)$ (the space itself) is attached trivially to $\varOmega_{(3,1)}$. On the other hand, we have
\[Sq^4\, \sigma_{(1,1)}=\sigma_{(2,2)}\,,\quad (Sq^2)^2\, \sigma_{(1,1)}=0\,;\]
in fact, the integral lift gives $(SQ^1\circ SQ^1)\,s_{(1,1)}=2 s_{(3,1)}\,$. However,
\[\calp^1\, \sigma_{(1,1)}=\sigma_{(3,1)}-\sigma_{(2,2)}\spa (\mbox{for } p=3),\]
which means that $\varOmega_{(3,1)}$ is attached non-trivially to $\varOmega_{(1,1)}$. 

The previous example confirms that the mod 2 Steenrod operations give only partial information about the attaching maps. However, they provide a good first approximation of the cell structure of a space. Then, one has to use more sophisticated methods to decide whether a cell which does not appear in any composition of Steenrod squares on a lower dimensional cell is attached trivially to the latter or not. For instance, one might start by considering Steenrod operations corresponding to primes $p>2$, or Adams operations in $K$-theory etc. It is also important to distinguish between unstable and stable attachment.

The approximation of the cell structure of a space provided by the Steenrod squares can be viewed as a poset structure on the set of cells. Indeed, we consider 
finger the transitive closure of the relation $\varOmega_\lambda\le\varOmega_\mu$ if and only if the Schubert class $\sigma_\mu$ appears with non-zero coefficient in some composition of Steenrod squares on $\sigma_\lambda$.  Clearly, this poset is a subposet of the set of cells ordered by inclusion. Let us also note that it is enough to consider only {\em admissible monomials} in the Steenrod squares (that is monomials $Sq^{i_1}\ldots Sq^{i_k}$ with $i_{j-1}\ge 2i_j$ for $1<j\le k$), since they form a basis of the mod 2 Steenrod algebra. Alternatively, one can use only the Steenrod squares $Sq^{2^k}$ for $k\ge 0$, since they generate the mod 2 Steenrod algebra. Carrying out the computations above for all the cells in $Gr_2(\bC^5)$, we can easily draw the Hasse diagram of the poset corresponding to this Grassmannian (see the figure below); the cells are represented by the Young diagrams corresponding to the partitions which index them.

\[
\begin{array}{c}
\mbox{\psfig{file=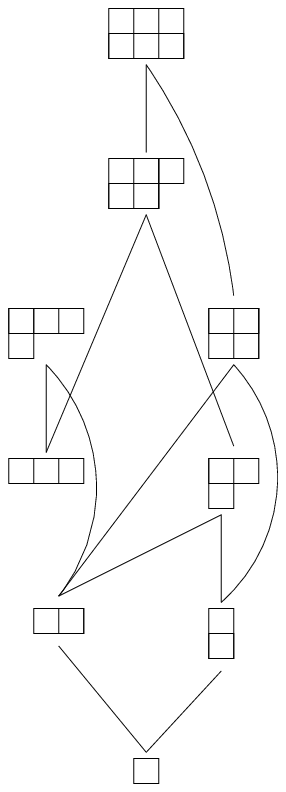}}
\end{array}
\]

\end{document}